\documentclass[11pt]{article}
\usepackage{mathrsfs}
\usepackage{amssymb}
\usepackage{amsmath}
\usepackage{amsbsy}
\usepackage{epsfig}
\usepackage{enumerate}
\usepackage{bm}
\usepackage{color}
\usepackage{booktabs}

\RequirePackage[OT1]{fontenc}

\RequirePackage[aos]{imsart}

\newtheorem{theorem}{Theorem}[section]
\newtheorem{proposition}{Proposition}[section]

\newtheorem{corollary}{Corollary}[section]
\newtheorem{remark}{\sc Remark}[section]

\def\proclaim#1{\par \smallskip\noindent {\bf #1}\bgroup\it\ }
\def\endproclaim{\egroup\par\smallskip}

\startlocaldefs
\endlocaldefs
\newbox\TempBox \newbox\TempBoxA

\def\pr{\textsf{P}} 
\def\ep{\textsf{E}} 

\def\bk#1{\bm #1} 

\def\underwiggle 1{
\ifmmode\setbox\TempBox=\hbox{$ 1$}\else\setbox\TempBox=\hbox{1}\fi
\setbox\TempBoxA=\hbox to \wd\TempBox{\hss\char'176\hss}
\rlap{\copy\TempBox}\smash{\lower9pt\hbox{\copy\TempBoxA}}}
\begin{document}

\begin{frontmatter}

\title{On the Theory  of Covariate-Adaptive  Designs}
\runtitle{Covariate-Adaptive Designs}

\begin{aug}

\author{\fnms{Feifang} \snm{HU}\thanksref{t2}\ead[label=e2]{fh6e@virginia.edu}}
and
\author{\fnms{Li-Xin} \snm{ZHANG}\thanksref{t1}\ead[label=e1]{stazlx@zju.edu.cn}}

\address{Feifang HU \\
DEPARTMENT OF STATISTICS \\
UNIVERSITY OF VIRGINA\\
HALSEY HALL, CHARLOTTESVILLE~~~~~~ \\
VIRGINIA 22904-4135, USA \\
\printead{e2}}

\address{Li-Xin ZHANG\\
DEPARTMENT OF MATHEMATICS$~~~~~~~~~~~~~~~$ \\
ZHEJIANG UNIVERSITY\\
ZHEDA ROAD, NO. 38~~~~~~ \\
HANG ZHOU, 310027~, P.R. China \\
\printead{e1}}

\affiliation{University of Virginia and Zhejiang University}
\thankstext{t2}{Research supported by grants DMS-0907297 and DMS-1209164 from the National Science Foundation (USA).}
\thankstext{t1}{Research supported by a grant from the National Natural Science Foundation of China (No. 11225104),  grants from Natural Science Foundation of Zhejiang Province (No. R6100119) and  Fundamental Research Funds  for the Central Universities. Corresponding author.}

\runauthor{F. Hu and L-X. Zhang }
\end{aug}

\begin{abstract}
Pocock and Simon's marginal procedure (Pocock and Simon, 1975) is often implemented for
balancing treatment allocation over influential covariates in clinical trials.
However, the theoretical properties of  Pocock and Simion's procedure have remained largely elusive for decades.
In this paper, we propose a general framework for covariate-adaptive designs and establish the corresponding theory under widely satisfied conditions. As a special case, we obtain the theoretical properties of Pocock and Simon's marginal procedure: the marginal imbalances and overall imbalance  are bounded in probability, but the within-stratum imbalances increase  with the rate of $\sqrt{n}$ as the sample size increases. The theoretical results provide new insights about balance properties of covariate-adaptive randomization procedures and open a door to study the theoretical properties of statistical inference for clinical trials based on covariate-adaptive randomization procedures.

\end{abstract}



\end{frontmatter}

\section{Introduction}
\label{s:intro}
\setcounter{equation}{0}

\begin{center}
\end{center}

\vspace{-0.3in}
It is well known that covariates  play an important role in clinical trials.
Clinical trialists  are often  concerned about unbalanced treatment arms
with respect to key covariates of interest.  In the literature,
covariate-adaptive randomization procedures (Rosenberger and Lachin, 2002) are sometimes employed to balance on important covariates. Pocock and Simon's marginal procedure (Pocock and Simon, 1975)
is one popularly used method in the literature. As pointed out in Taves (2010), there are over 500 clinical trials which implemented Pocock and Simon's marginal procedure to balance important categorical covariates from 1989 to 2008. Simulation studies (Weir and Lees, 2003; Toorawa, Adena {\em et al.}, 2009; Kundt, 2009)
found that Pocock and Simon's marginal procedure indeed reduces marginal
imbalances as well as the overall imbalance.
But the performance within  strata is not as satisfactory (Signorini, Leung \emph{et al.}, 1993; Kundt, 2009). However,  all studies of Pocock and Simon's procedure are merely carried out by simulations.
There is ``no theoretical justification that the procedure even works as intended" (Rosenberger
and Sverdlov, 2008).

Over the past several decades,  scientists have identified many new biomarkers (Ashley, Butte, Matthew {\em et al.}, 2010;  Lipkin, Chao  {\em et al.} , 2010;
McIlroy,  McCartan {\em et al.},  2010; etc.) that may link with certain diseases in the fields of translational research
(genomics, proteomics, and metabolomics).  Based on these biomarkers, we would like develop
personalized medicine algorithms that help patients to receive better treatment regimens based on
their individual characteristics (which could be biomarkers or other covariates).	
To design a superior and efficient clinical study for personalized medicine, one should incorporate information on important biomarkers (Hu, 2012). Therefore, balancing treatment allocation for influential covariates has become increasingly important in today's clinical trials.
As pointed out in Hu (2012), classical covariate-adaptive
designs have several drawbacks to incorporate many important biomarkers.
Recently,  Hu and Hu (2012) developed a class of covariate-adaptive biased coin randomization procedure and studied its theoretical properties. However,  condition (C) of their Theorem 3.2 is very strict and hard to  verify in practice. More importantly, Hu and Hu's theoretical results do not apply to Pocock and Simon's marginal procedure.

In this paper, we  establish a general theoretical foundation for covariate-adaptive randomization procedures under widely satisfied conditions. In particular, we have theoretically proved that under Pocock and Simon's marginal procedure, the marginal imbalances and overall imbalance  are bounded in probability, but  all the within-stratum imbalances increase  with the rate of $\sqrt{n}$ as the sample size increases. The theory provides some new insights about theoretical properties of covariate-adaptive randomization procedures. In particular, the theory provides critical conditions for general covariate-adaptive randomization procedures to achieve good within-stratum balance and good marginal balance. As discussed in the concluding remarks (Section 5), our theoretical results also open a door to study the
theoretical behavior of inferential methods (estimation, hypothesis testing, etc.) of clinical trials based on covariate-adaptive randomization procedures.

We first propose a general family of covariate-adaptive randomization procedures which includes some new designs and many  existing designs as special cases: stratified randomization, Pocock and Simon's (1975) marginal procedure, Hu and Hu's (2012) procedures, Efron's (1971) biased coin design and Baldi Antognini  and Giovagnoli's  (2004) adjustable biased coin design, etc. To study the theoretical properties under this general framework,
the main difficulties include (i) the correlation structure of within-stratum imbalances; (ii) the relationship among within-stratum and marginal imbalances under Pocock and Simon's type procedures; and (iii) the discreteness of the allocation function. In the literature, Taylor's expansion and martingale approximation are two common techniques to study the theoretical properties of adaptive designs with a continuous allocation function (Hu and Zhang, 2004; Hu and Rosenberger, 2006;  Zhang, Hu and Cheung, 2006). Under Pocock and Simon's type procedures, to overcome the complex relationship among within-stratum and marginal imbalances, we have to approximate these imbalances using martingales by solving Poisson's equations. To deal with the discreteness of the allocation function, we use the technique of ``drift conditions'' (Meyn and Tweedie 1993), which was developed for Markov chains on general state spaces.

The paper is organized as following.
The general framework of the randomization procedure is described
in Section \ref{s:procedure} and the theoretical results  are given in
Section \ref{s:theory}.  We consider multi-arm clinical trials in Section 4.
Some concluding remarks are in Section 5.
The proofs of the theorems can be found in Section \ref{s:proofs}.

\section{The Framework for Covariate-Adaptive Designs}
\label{s:procedure}
\setcounter{equation}{0}
 \begin{center}
\end{center}

\vspace{-0.3in}

We consider the same  setting  as that of Pocock and Simon (1975) and only focus on two treatment groups $1$ and $2$.
Consider $I$ covariates and $m_{i}$ levels for the $i$th covariate, resulting in $m=\prod_{i=1}^{I}m_{i}$ strata.
Let $T_{j}$ be the assignment of the $j$th patient, $j=1,\ldots,n$,
i.e., $T_{j}=1$ for treatment 1 and $T_{j}=0$ for treatment 2.
Let $Z_{j}$ indicate the covariate profile of that patient,
i.e., $Z_{j}=(k_{1},\ldots,k_{I})$ if his or her $i$th covariate is at level $k_{i}$, $1\leq i \leq I $ and $1\leq k_{i}\leq m_{i}$.
For convenience, we use $(k_{1},\ldots,k_{I})$ to denote the \emph{stratum} formed by patients who possess the same covariate profile $(k_{1},\ldots,k_{I})$, and use $(i;k_{i})$ to denote the \emph{margin} formed by patients whose $i$th covariate is at level $k_{i}$.

The  procedure is defined as follows:
\begin{enumerate}
\item[1)] The first patient is assigned to treatment 1 with probability 1/2.

\item[2)] Suppose $(n-1)$ patients have been assigned to  treatments ($n>1$) and
the $n$th patient falls within  stratum $(k_{1}^{*},\ldots,k_{I}^{*})$.

\item[3)] For the first $(n-1)$ patients,
     \begin{itemize}
     \item[-]let $D_{n-1}$ be the difference between the numbers of patients in treatment group 1 and 2, i.e., the number in group 1 minus the number in group 2;
     \item[-]similarly, let $D_{n-1}(i;k_{i}^{*})$ and $D_{n-1}(k_{1}^{*},\ldots,k_{I}^{*})$ be the differences between the numbers of patients in the two treatment groups on the margin $(i;k_{i}^{*})$, and within the stratum $(k_{1}^{*},\ldots, k_{I}^{*})$, respectively;
     \item[-] these differences can be positive, negative or zero, and each one is used to measure the \emph{imbalance} at the corresponding level (overall, marginal, or within-stratum).
     \end{itemize}

\item[4)] If the $n$th patient were assigned to treatment $1$, then $D^{(1)}_{n}=D_{n-1}+1$ would
be the ``potential'' overall difference in the two groups; similarly,
\begin{displaymath}
D^{(1)}_{n}(i;k_{i}^{*})= D_{n-1}(i;k_{i}^{*})+1
\end{displaymath}
and
\begin{displaymath}
D^{(1)}_{n}(k_{1}^{*},\ldots,k_{I}^{*})=D_{n-1}(k_{1}^{*},\ldots,k_{I}^{*})+1
\end{displaymath}
would be the potential differences on  margin $(i;k_{i}^{*})$ and within  stratum  $(k_{1}^{*},\ldots,k_{I}^{*})$, respectively.

\item[5)] Define an imbalance measure $Imb_{n}^{(1)}$ by
\begin{displaymath}
Imb_{n}^{(1)}=w_{o}[D^{(1)}_{n}]^2+\sum_{i=1}^{I}w_{m,i}[D^{(1)}_{n}(i;k_{i}^{*})]^2
+w_{s}[D^{(1)}_{n}(k_{1}^{*},\ldots,k_{I}^{*})]^2,
\end{displaymath}
      which is the weighted imbalance that would be caused if the $n$th patient
      were assigned to treatment 1. $w_{o}$, $w_{m,i}$ $(i=1,\ldots,I)$ and $w_{s}$ are nonnegative
      weights placed on overall, within a covariate margin and within a stratum
      cell, respectively. Without loss of generality we can assume
\begin{displaymath}
w_{o}+w_{s}+\sum_{i=1}^{I}w_{m,i}=1.
\end{displaymath}
\item[6)] In the same manner we can define $Imb_{n}^{(2)}$, the weighted imbalance that would
be caused if the $n$th patient were assigned to treatment 2. In this case, the three types
of potential differences are the existing ones minus 1, instead of plus 1.

\item[7)]  Conditional on the assignments of the first $(n-1)$ patients as well as the covariates' profiles of the first $n$ patients, assign the $n$th patient to treatment 1 with probability
\begin{align}\label{eqallocationP}
&P(T_{n}=1|\bm{Z}_{n-1},Z_n=(k_1^{\ast},\ldots, k_I^{\ast}),\bm{T}_{n-1})\nonumber\\
& \quad = g\left(Imb_{n}^{(1)}-Imb_{n}^{(2)}\right)
\end{align}
  where $n>1$,  $\bm{Z}_{n-1}=(Z_{1},\ldots,Z_{n-1})$,
  $\bm{T}_{n-1}=(T_{1},\ldots,T_{n-1})$,   $g(x)$ is a   real  function with $0< g(x)< 1$, $g(-x)=1-g(x)$,  $$ g(x)\le 0.5   \text{ when } x\ge 0, \; \text{ and }\; \limsup_{x\to +\infty} g(x)<0.5.
  $$
  \end{enumerate}

Using the basic equation $(x+1)^2-(x-1)^2=4x$, the critical quantity $Imb_{n}^{(1)}-Imb_{n}^{(2)}$ in Step $7)$  can be simplified as
\begin{align}
&Imb_{n}^{(1)}-Imb_{n}^{(2)}\nonumber\\
=&4\left\{w_{o}D_{n-1}+\sum_{i=1}^{I}w_{m,i}D_{n-1}(i;k_{i}^{*})+w_{s}D_{n-1}(k_{1}^{*},\ldots,k_{I}^{*})\right\}\nonumber\\
:=&4\cdot\Lambda_{n-1}(k_1^*,\ldots,k_I^*) \label{simplification}
\end{align}
Therefore, the allocation probability $g\left(Imb_{n}^{(1)}-Imb_{n}^{(2)}\right)$ is determined by the value of $\Lambda_{n-1}(k_1^*,\ldots,k_I^*)$, which is a weighted average of current imbalances at different levels. In the literature different views have been given as to the selection of the allocation probability function $g(\cdot)$. Efron (1971), Pocock and Simon (1975), Hu and Hu (2012) suggested
\begin{equation} \label{eqBCDfunction} g(x)=\begin{cases} q, & \text{ if } x>0,\\
\frac{1}{2},&  \text{ if } x=0,\\
p, & \text{ if } x<0,
\end{cases}
\end{equation}
where $p>1/2$ and $q+p=1$. In general, we can define $g(x)$ to be either a continuous function or a discrete function.

\begin{remark} Instead of using the biased coin function (\ref{eqBCDfunction}) (Efron, 1971; Pocock and Simon, 1975; etc.), we use a general allocation function  which  is defined as a decreasing function of the weighted average imbalances. When the covariates are not considered,  Baldi Antognini  and Giovagnoli  (2004) suggested a heavy tail function which can reduce both the allocation  bias and selection bias  of Efron's biased coin design. We hope that  the general framework is flexible enough to define applicable randomization procedures with good properties. The theoretical results will be established under widely satisfied  conditions  so that they can apply for  all cases. In practice, one may use Efron's biased coin function (\ref{eqBCDfunction}) with $p \in [0.75, 0.95]$ as discussed and suggested in the literature (Hu and Hu, 2012). However, selection bias could be a concern with a large $p$ in Efron's biased coin function (\ref{eqBCDfunction}). To reduce selection bias and allocation bias, one may use the heavy tail function suggested by Baldi Antognini  and Giovagnoli  (2004).
\end{remark}

\section{Theoretical Properties}
\label{s:theory}
\setcounter{equation}{0}

\begin{center}
\end{center}

\vspace{-0.3in}
We now investigate the asymptotic properties of the design. For the first $n$ patients, we know that $D_{n}(k_{1},\ldots,k_{I})$
is the true difference between the two treatment arms   within stratum $(k_{1},\ldots,k_{I})$. Let
$$
\bm {D}_{n}=\left[D_{n}(k_{1},\ldots,k_{I})\right]_{1\leq k_{1}\leq m_{1},\ldots, 1\leq k_{I}\leq m_{I}}
$$
be an array of dimension $m_{1}\times\ldots \times m_{I}$ which stores the current assignment
differences in all strata and therefore stores the current imbalances.
Also, assume that the covariates $Z_{1},Z_{2},\ldots$ are independently and identically
distributed. Since $Z_{n}=(k_{1},\ldots, k_{I})$ can take
 $m=\prod_{i=1}^{I} m_i$  different values, it  in fact follows an $m$-dimension
multinomial distribution with parameter $\bm p=(p(k_{1},\ldots, k_{I}))$, each element being the probability that a patient falls within the corresponding stratum. Obviously,  $p(k_{1},\ldots, k_{I}) \geq 0 $ and $\sum_{k_{1},\ldots, k_{I}} p(k_{1},\ldots, k_{I})=1$.
Notice
$$ D_n(k_1,\ldots,k_I)=D_{n-1}(k_1,\ldots,k_I)+2\big(T_n-\frac{1}{2}\big)\mathbb{I}\{Z_n=(k_1,\ldots,k_I)\}. $$
It is easily seen that
$$ \pr(D_n(k_1,\ldots,k_I)=D_1(k_1,\ldots,k_I) \; \forall n)=1 \;\; \text{ if } p(k_1,\ldots, k_I)=0. $$
We can ignore those strata with $p(k_1,\ldots, k_I)=0$. Hence without loss of generality, we assume $p(k_1,\ldots, k_I)>0$ for all $(k_1,\ldots, k_I)$.

Our purpose is to the study the properties of $\bm D_n$. Besides $\bm D_n$, we will also consider the weighted average of the imbalances $\bm\Lambda_{n-1}$ as in (2.1). Let
 \begin{align*}
\Lambda_{n}(k_1,\ldots,k_I)=& w_{o}D_{n}+\sum_{i=1}^{I}w_{m,i}D_{n}(i;k_{i})+w_{s}D_{n}(k_{1},\ldots,k_{I}),
\\
 \bm \Lambda_{n}=& \left[\Lambda_{n}(k_{1},\ldots,k_{I})\right]_{1\leq k_{1}\leq m_{1},\ldots, 1\leq k_{I}\leq m_{I}}.
\end{align*}
The allocation probability (\ref{eqallocationP}) of the $n$-th patient  is a function of $\bm \Lambda_{n-1}$. We will find later that $\bm \Lambda_n$ plays a very important role for investigating the properties of $\bm D_n$.    It is obvious that $\bm \Lambda_n=\bm L(\bm D_n):\bm D_n \to \bm \Lambda_n$ is a linear transform of $\bm D_n$.
The following proposition gives the relation between $\bm D_n$ and $\bm\Lambda_n$ and tells us that both $(\bm D_n)_{n\ge 1}$ and $(\bm\Lambda_n)_{n\ge 1}$ are Markov chains.
\begin{proposition}\label{proposition1} (i) If $w_s>0$, then $\bm \Lambda_n=\bm L(\bm D_n)$ is a one to one linear map; If $w_s+w_{m,i}>0$, then each $D_n(i;k_i)=D_{i;k_i}(\bm\Lambda_n)$ is a linear transform of $\bm \Lambda_n$; For any case, $D_n=D(\bm\Lambda_n)$ is a linear transform of $\bm \Lambda_n$;

(ii) $(\bm{D}_{n})_{n\geq 1}$ is an irreducible  Markov chain on the space $\mathbb{Z}^{m}$ with period 2 and with the property that $(-\bm D_n)_{n\geq 1}$ and $(\bm D_n)_{n\geq 1}$ have the same transition probabilities;

(iii)  $(\bm \Lambda_n)_{n\geq 1}$ is an irreducible Markov chain on the space $\bm L(\mathbb{Z}^{m})$ with period 2 and with the property that $(-\bm \Lambda_n)_{n\geq 1}$ and $(\bm \Lambda_n)_{n\geq 1}$ have the same transition probabilities.

\end{proposition}

Now we give our main results.
\begin{theorem}\label{theorem1}
Consider $I$ covariates and $m_{i}$ levels for the $i$th covariate, where $I\geq1$, $1\leq i \leq I$, and $m_{i}>1$. $w_{o}$, $w_{s}$, and $w_{m,i}$, $i=1,\ldots,I$, are nonnegative with $w_{o}+\sum_{i=1}^{I}w_{m,i}+w_{s}=1$.
Then $(\bm\Lambda_n)_{n\ge 1}$ is a positive recurrent Markov chain with period 2 on $\bm L(\mathbb{Z}^{m})$ and $\ep \|\bm \Lambda_{n}\|^r=O(1)$ for any $r>0$. In particular,
\begin{description}
  \item[\rm (i)] If $w_s>0$, then $(\bm D_n)_{n\ge 1}$ is a positive recurrent Markov chain with period 2 on $\mathbb{Z}^{m}$, and $\ep \| \bm{D}_{n}\|^r=O(1)$ for any $r>0$;
  \item[\rm (ii)]  If $w_s+w_{m,i}>0$, then $D_n(i;k_i)=O(1)$ in probability and $\ep\left|D_n(i;k_i)\right|^r=O(1)$ for any $r>0$; Further, if $w_s=0$, then the collection of all marginal imbalances  $\big((D_n(i;k_i):w_{m,i}\ne 0, k_i=1,\ldots,m_i,i=1,\ldots, I)\big)_{n\ge 1}$ is  a positive recurrent Markov chain;
  \item[\rm (iii)]  For any case $D_n=O(1)$ in probability and $\ep|D_n|^r=O(1)$ for any $r>0$; Further, if $w_s=w_{m,1}=\ldots=w_{m,I}=0$, then  $(D_n)_{n\ge 1}$ a positive recurrent Markov chain.
\end{description}
\end{theorem}

\begin{remark}
Recently Hu and Hu (2012) obtained theoretical result (i) under very strict condition of the weights $w_s$ and $w_{m,i}$ when $g(x)$ is defined  in (\ref{eqBCDfunction}). The Condition (C)  in their Theorem 3.2   is very restrictive and usually not satisfied in practice. When the number of strata is large, their Condition (C)  can be satisfied only when $w_s$ is very close to $1$ and the design reduces closely to  stratified randomization. Their results do not apply to Pocock and Simon's (1975) design (where $w_s=0$) and the design with equal weights $w_0,w_{m,i},w_s$. Our Theorem \ref{theorem1} eliminates Hu and Hu (2012)'s  condition (C) so that it applies to most covariate-adaptive randomization procedures.
\end{remark}

The next theorem tell us that the within-stratum imbalances $|D_n(k_1,\ldots, k_I)|$ either are bounded in probability or increase with rate $\sqrt{n}$ as the sample increases.
\begin{theorem}\label{theorem2} Under the conditions in Theorem \ref{theorem1},
\begin{description}
  \item[\rm (iv)] There exist non-negative constants $\sigma(k_1,\ldots,k_I)$ such that
  \begin{equation}\label{varainceofD1}
  \ep D_n^2(k_1,\ldots,k_I)=n \sigma^2(k_1,\ldots,k_I) +O(\sqrt{n}\,\sigma(k_1,\ldots,k_I)),
    \end{equation}
  \begin{equation}\label{eqCLTofD}
  \frac{ D_n(k_1,\ldots,k_I)}{\sqrt{n}}\overset{D}\to N\big(0,\sigma^2(k_1,\ldots,k_I)\big)
  \end{equation}
  and
 \begin{equation}\label{eqth1.3}
  \lim_{n\to \infty} \ep\left|\frac{ D_n(k_1,\ldots,k_I)}{\sqrt{n}}\right|^r=\sigma^r(k_1,\ldots,k_I) \ep\left|N(0,1)\right|^r
  \end{equation}
  for all stratum $(k_1,\ldots,k_I)$s and $r>0$, where $N(0,1)$ is a standard normal random variable;
  \item[\rm (v)]  For any fixed stratum $(k_1,\ldots,k_I)$, if $ D_n(k_1,\ldots,k_I)=o(\sqrt{n})$ in probability, then $D_n(k_1,\ldots,k_I)=O(1)$ in probability;
  \item[\rm (vi)] If $ D_n(k_1,\ldots,k_I)=o(\sqrt{n})$ in probability for one stratum $(k_1,\ldots,k_I)$, then $w_s\ne 0$. In other words, if $w_s=0$, then for all stratum  $(k_1,\ldots,k_I)$
      $$ \lim_{n\to \infty} \frac{\ep D_n^2(k_1,\ldots,k_I)}{n}=\sigma^2(k_1,\ldots,k_I)>0. $$
\end{description}
\end{theorem}

The main conclusion of Theorems \ref{theorem1} and \ref{theorem2} can be summarized in the following corollary which indicates that the condition $w_s>0$ is critical to ensure that $(\bm{D}_{n})_{n\geq 1}$ is positive recurrent.
\begin{corollary} The following statements are equivalent:
\begin{description}
  \item[\rm (1)]  $(\bm D_n)_{n\ge 1}$ is a positive recurrent Markov chain;
  \item[\rm (2)] $\bm D_n=O(1)$  in probability;
  \item[\rm (3)]  $\ep\|\bm D_n\|^r=O(1)$  for all $r>0$;
  \item[\rm (4)] $D_n(k_1,\ldots,k_I)=o(\sqrt{n})$  in probability for at least one stratum $(k_1,\ldots,k_I)$;
  \item[\rm (5)] $w_s>0$.
\end{description}
\end{corollary}

The next theorem tells us that the marginal procedures will not provide good balance with respect to the margin if the margin is not considered in the imbalance measure for defining the allocation probability (\ref{eqallocationP}).

\begin{theorem}\label{theorem3} Suppose the conditions in Theorem \ref{theorem1} are satisfied. If $w_s+w_{m,i}=0$, then
$$ \lim_{n\to \infty} \frac{\ep[D_n^2(i;k_i)]}{n}>0\;\; \text{ for all } k_i=1,\ldots, m_i. $$
\end{theorem}

\begin{remark} By Theorem \ref{theorem1}, \ref{theorem2} and Theorem \ref{theorem3}, the conditions $w_s>0$ and $w_s+w_{m,i}>0$ ($i=1,...,I$) are critical to ensure that the within stratum $D_n(k_1,...,k_I)=O(1)$ and the marginal imbalances $D_n(i;k_i)=O(1)$ in probability respectively. However, we have not discussed the selection of these weights in practice.
Here are some suggestions based on the results of this paper: (i) Always choose $w_s>0$. (ii) When the sample size is relatively large and the total number of strata is relatively small, there are enough patients in each strata.
In these cases, balance within each strata is important and $w_s$ plays an important role, we may choice a relatively large $w_s$. For example, we may use $w_s=1/2$ under these situations.
(iii) When the number of covariates ($I$) are increasing and the number of strata is relatively large, we may select weights according to the number of covariates ($I$) and the important of each covariate. For example, we may select $w_s=w_{m,i}=(I+1)^{-1}$ or $w_s=(I+1)^{-1}$ and $w_{m,i}$ according the important of $i$th covariate ($i=1,...,I$). Some simulation studies can be found in Hu and Hu (2012).
\end{remark}

It is an interesting observation from Theorem \ref{theorem2} that if one of the asymptotic variances $\sigma^2(k_1,\ldots, k_I)$ is positive, then all of them are positive, while, if they are zeros, then $\bm D_n$ is bounded in probability. The later will happen only in the case of $w_s>0$.
When $w_s=0$, the design reduces to the marginal procedure which includes Pocock and Simon's (1975) design as a special case. Based on Theorem \ref{theorem1} (ii), Theorem \ref{theorem2} (iv) and (vi), and Theorem \ref{theorem3}, we have the following asymptotic properties of Pocock and Simon's procedure.

\begin{corollary}
For Pocock and Simon's marginal procedure ($w_s=0$), we have the following results:
\begin{description}
\item[\rm (a)] All within-stratum imbalances  increase  with the rate $\sqrt{n}$ as the sample size increases. Also $D_n(k_1,\ldots,k_I)/\sqrt{n}$ is asymptotically normal distributed with a positive variance $\sigma^2(k_1,\ldots,k_I)$.
\item[\rm (b)] When $w_{m,i}>0$, then the corresponding marginal imbalance (the $i$-th covariate) and the overall imbalance are bounded in probability, that is, $D_n(i;k_i)=O(1)$ and $D_n=O(1)$ in probability; Further,  the collection of all marginal imbalances  $\big((D_n(i;k_i):w_{m,i}\ne 0, k_i=1,\ldots,m_i,i=1,\ldots, I)\big)_{n\ge 1}$ is  a positive recurrent Markov chain with period 2.
\item[\rm (c)] When $ w_{m,i}=0$,  then the corresponding marginal imbalance increase  with the rate $\sqrt{n}$, that is, $D_n(i;k_i)=O(\sqrt{n})$ in probability.
\end{description}
\end{corollary}

As in Hu and Hu (2012), to prove the Theorem \ref{theorem1}, we will use the technique of ``drift conditions'' (Meyn and Tweedie 1993), which was developed for Markov chains on general state spaces.  In stead of considering $(\bm D_n)$ directly as in Hu and Hu (2012), we have to consider $(\bm \Lambda_n)$ in this paper.  In order to prove the positive recurrence of $(\bm \Lambda_n)_{n\geq 1}$ we need to find a test function
   $V: \bm L(\mathbb{Z}^{m})\rightarrow \mathbb{R}^{+}$, a bounded test set $\mathscr{C}$ on $\bm L(\mathbb{Z}^{m})$,
   and a positive constants $b$ such that
\begin{equation}
\bigtriangleup_{\lambda} V (\bm \Lambda):= \sum_{\bm \Lambda'\in\bm L(\mathbb{Z}^{m})}P_{\lambda}(\bm \Lambda,\bm \Lambda')V(\bm \Lambda')-V(\bm \Lambda)\label{deltavd}
\end{equation}
satisfies the following  condition:
\begin{align}
\bigtriangleup_{\lambda} V (\bm \Lambda)\leq -1+b\mathbb{I}_{\bm \Lambda\in \mathscr{C}},\label{drift1}
\end{align}
where $P_{\lambda}(\bm \Lambda,\bm \Lambda')$ is the transition probability from $\bm \Lambda$ to $\bm \Lambda'$
on the state space $\bm L(\mathbb{Z}^{m})$ of the chain $(\bm \Lambda_n)_{n\geq 1}$, and $\mathbb{I}_{\bm\Lambda\in \mathscr{C}}$  is  a function with value $1$ if $\bm\Lambda$ is in $\mathscr{C}$, and zero if not. $V$ is often a norm-like function
on $\bm L(\mathbb{Z}^{m})$. For considering the convergence of moments of the Markov chain, we will also find the drift condition of $\bigtriangleup_{\lambda} V^r(\bm\Lambda)$.  The test function $V$ is the key component  in the proofs. We have to choose a good $V$ such that it is norm-like function and the drift $\bigtriangleup_{\lambda} V$ is also very close to the norm of $\bm\Lambda$, so that the drift condition is satisfied without any additional condition on the weights $w_o$, $w_s$, and $ w_{m,i}$s.

When $w_s=0$ (Pocock and Simon's marginal procedure), the within-stratum imbalance $D_n(k_1,\ldots, k_I)$ is not considered in the allocation procedure. We need to introduce a new technique (Poisson's equation) to deal with the complicated structure of the within-stratum imbalances and marginal imbalances. In fact,
we will approximate $D_n(\bm k)$ as a martingale plus a function of $\bm \Lambda_n$ by solving   Poisson's equation in the proof of Theorem \ref{theorem2}.  We will prove that this martingale is a constant when the asymptotic variance $\sigma(k_1,\ldots,k_I)$ is zero and so that $D_n(k_1,\ldots, k_I)$ is a function of $\bm \Lambda_n$,  which is a contradiction when $w_s=0$.  All the proofs    will stated in the Section 6.

\section{An extension to the multi-arm clinical trials}
\label{s:multi}
\setcounter{equation}{0}
\begin{center}
\end{center}

\vspace{-0.3in}
In some clinical trials,  one would like to compare three or more treatments (Pocock and Simon, 1975; Tymofyeyev, Rosenberger and Hu, 2007; Hu, 2012; etc.).
In this section we consider clinical trials with $T$ ($T>2$) treatments. Let $T_{j}$ be the assignment of the $j$th patient, $j=1,\ldots,n$,
i.e., $T_{j}=t$ for treatment $t$. Under the same covariate structure of Section 2,
the allocation procedure is defined as follows:
\begin{enumerate}
\item[1)] The first patient is assigned to treatment $t$ with probability $1/T$.

\item[2)] Suppose $(n-1)$ patients have been assigned to  treatments ($n>1$) and
the $n$-th patient falls within  stratum $(k_{1}^{*},\ldots,k_{I}^{*})$.

\item[3)] For the first $(n-1)$ patients,
  let $N_{n-1,t}$ be the  number of patients in treatment group t.  Further let $N_{n-1,t}(i;k_i^{\ast})$ and $N_{n-1,t}(k_1^{\ast},\ldots,k_I^{\ast})$ be the numbers in treatment group $t$ on the margin $(i;k_i^{\ast})$ and within the stratum $(k_1^{\ast},\ldots,k_I^{\ast})$, respectively. We denote
 $$N_{n-1}^{Ave}=\frac{1}{T}\sum_{t=1}^T N_{n-1,t}, ~N_{n-1}^{Ave}(i;k_i^{\ast})= \frac{1}{T}\sum_{t=1}^TN_{n-1,t}(i;k_i^{\ast}),$$
 $$\mbox{ and }N_{n-1}^{Ave}(k_1^{\ast},\ldots,k_I^{\ast})=\frac{1}{T}\sum_{t=1}^T N_{n-1,t}(k_1^{\ast},\ldots,k_I^{\ast})$$ be  the corresponding average numbers over treatments.
         Define the differences
         \begin{align*}
         & D_{n-1,t}= N_{n-1,t}-N_{n-1}^{Ave},\\
         & D_{n-1,t}(i,k_i^{\ast})=N_{n-1,t}(i,k_i^{\ast})-N_{n-1}^{Ave}(i,k_i^{\ast}), \mbox{ and } \\
         & D_{n-1,t}(k_1^{\ast},\ldots,k_I^{\ast})=N_{n-1,t}(k_1^{\ast},\ldots,k_I^{\ast})-N_{n-1}^{Ave}(k_1^{\ast},\ldots,k_I^{\ast}).
         \end{align*}
These differences  are used to measure the overall imbalance, the imbalance  on the margin $(i;k_i^{\ast})$ and  the imbalance  within stratum $(k_1^{\ast},\ldots,k_I^{\ast})$, respectively, for each treatment $t$.

\item[4)] If the $n$-th patient is assigned to treatment $t$, then $N_{n-1,t}$, $N_{n-1,t}(i;k_i^{\ast})$ and $N_{n-1,t}(k_1^{\ast}, \ldots, k_I^{\ast})$ will increase $1$, and others remain unchanged. So the ``potential'' imbalance  at the corresponding level (overall, marginal, and within-stratum) is
\begin{align*}
& D^{(t)}_{n,h}=D_{n-1,h}+\mathbb I\{h=t\}-\frac{1}{T},\\
& D^{(t)}_{n,h}(i;k_{i}^{*})= D_{n-1,h}(i;k_{i}^{*})+\mathbb I\{h=t\}-\frac{1}{T}, \mbox{ and }\\
& D^{(t)}_{n,h}(k_{1}^{*},\ldots,k_{I}^{*})=D_{n-1,h}(k_{1}^{*},\ldots,k_{I}^{*})+\mathbb I\{h=t\}-\frac{1}{T},
\end{align*}
for $h=1,\ldots,T$.

\item[5)] Define an imbalance measure $Imb_{n,t}$ by
\begin{align*}
Imb_{n,t}=\sum_{h=1}^T\big\{ &w_{o}[D^{(t)}_{n,h}]^2+\sum_{i=1}^{I}w_{m,i}[D^{(t)}_{n,h}(i;k_{i}^{*})]^2 \\&
+w_{s}[D^{(t)}_{n,h}(k_{1}^{*},\ldots,k_{I}^{*})]^2\big\},
\end{align*}
      which is the weighted imbalance that would be caused if the $n$-th patient
      were assigned to treatment $t$. Here $w_{o}$, $w_{m,i}$ ($i=1,\ldots,I$) and $w_{s}$ are nonnegative
      weights with
$
w_{o}+w_{s}+\sum_{i=1}^{I}w_{m,i}=1.
$

\item[6)] Having the imbalance measure, we define the allocation probabilities  in the same way as Pocock and Simon (1975).
       One can rank the treatments according to the values of $Imb_{n,t}$, $t=1,\ldots T$, in a non-decreasing order so that
     $$ Imb_{n,(1)}\le Imb_{n,(2)}\le \ldots \le Imb_{n,(T)}. $$
    In the case of ties a random ordering can be determined.  The assigned treatment $T_n$ of  the $n$-th patient can be determined from the following set of probabilities
    \begin{align}\label{eqallocationPmul}
 P(T_{n}=(t)|\bm{Z}_{n-1},Z_n=(k_1^{\ast},\ldots, k_I^{\ast}),\bm{T}_{n-1})
  = p_t,
\end{align}
  where $n>1$,  $\bm{Z}_{n-1}=(Z_{1},\ldots,Z_{n-1})$,
  $\bm{T}_{n-1}=(T_{1},\ldots,T_{n-1})$, and    $p_1\ge p_2\ge \ldots \ge p_T$ are $T$ ordered nonnegative fixed constants with $\sum_t p_t=1$ and $p_1>p_T>0$.
  \end{enumerate}

  When  $w_{o}=w_{s}=0$, i.e., only the marginal imbalances are considered, the proposed design reduces to Pocock and Simon's (1975) marginal method.

  Let
$$
\bm {D}_{n}=\left[D_{n,t}(k_{1},\ldots,k_{I})\right]_{1\leq t\leq T, 1\leq k_{1}\leq m_{1},\ldots, 1\leq k_{I}\leq m_{I}},
$$
$$
\bm {\Lambda}_{n}=\left[\Lambda_{n,t}(k_{1},\ldots,k_{I})\right]_{1\leq t\leq T, 1\leq k_{1}\leq m_{1},\ldots, 1\leq k_{I}\leq m_{I}},
$$
where $\Lambda_{n,t}(k_1,\ldots,k_I)= w_{o}D_{n,t}+\sum_{i=1}^{I}w_{m,i}D_{n,t}(i;k_{i})+w_{s}D_{n,t}(k_{1},\ldots,k_{I})$. The following theorem is the main result for multi-treatment case.

\begin{theorem}\label{theorem1mul} Consider $I$ covariates and $m_{i}$ levels for the $i$th covariate, where $I\geq1$, $1\leq i \leq I$, and $m_{i}>1$. $w_{o}$, $w_{s}$, and $w_{m,i}$, $i=1,\ldots,I$, are nonnegative weights with $w_{o}+\sum_{i=1}^{I}w_{m,i}+w_{s}=1$. Assume that $p_1\ge p_2\ge \ldots \ge p_T$ are nonnegative constants with $\sum_t p_t=1$ and $p_1>p_T>0$.

Then $(\bm\Lambda_n)_{n\ge 1}$ is a positive recurrent Markov chain with period $T$   and $\ep \|\bm \Lambda_{n}\|^r=O(1)$ for any $r>0$. In particular,
\begin{description}
  \item[\rm (i)] If $w_s>0$, then $(\bm D_n)_{n\ge 1}$ is a positive recurrent Markov chain with period $T$, and $\ep \| \bm{D}_{n}\|^r=O(1)$ for any $r>0$;
  \item[\rm (ii)]  If $w_s+w_{m,i}>0$, then $D_{n,t}(i;k_i)=O(1)$ in probability and $\ep\left|D_{n,t}(i;k_i)\right|^r=O(1)$  for any $r>0$ and  $t=1,\ldots,T$; Further, if $w_s=0$, then the collection of all marginal imbalances  $\big(D_{n,t}(i;k_i):w_{m,i}\ne 0, t=1,\ldots, T, k_i=1,\ldots,m_i,i=1,\ldots, I)\big)_{n\ge 1}$ is  a positive recurrent Markov with period $T$ ;
  \item[\rm (iii)]  For any case $D_{n,t}=O(1)$ in probability and $\ep|D_{n,t}|^r=O(1)$ for any $r>0$ and $t=1,\ldots,T$; Further, if $w_s=w_{m,1}=\ldots=w_{m,I}=0$, then  $(D_{n,t};t=1,\ldots,T)_{n\ge 1}$ a positive recurrent Markov chain with period $T$.
\end{description}
Further,
\begin{description}
  \item[\rm (iv)]  if $w_s=0$, then for any stratum  $(k_1,\ldots,k_I)$ and treatment $t$,
      $$ \lim_{n\to \infty} \frac{\ep D_{n,t}^2(k_1,\ldots,k_I)}{n}>0; $$
 \item[\rm (v)]     if $w_s+w_{m,i}=0$, then
$$ \lim_{n\to \infty} \frac{\ep[D_{n,t}^2(i;k_i)]}{n}>0\;\; \text{ for all } k_i=1,\ldots, m_i, t=1,\ldots,T. $$
\end{description}
\end{theorem}

When $w_s=0$, the design reduces to the marginal procedure which includes Pocock and Simon's (1975) design as a special case. Based on Theorem \ref{theorem1mul}, we have the following asymptotical properties of Pocock and Simon's procedure.

\begin{corollary}
For Pocock and Simon's marginal procedure ($w_s=0$ and $w_0=0$), we have the following results:
\begin{description}
\item[\rm (a)] All within-stratum imbalances  increase  with the rate $\sqrt{n}$ as the sample size increases.
\item[\rm (b)] When $w_{m,i}>0$, then the corresponding marginal imbalance (the $i$-th covariate) and the overall imbalance are bounded in probability, that is, $D_{n,t}(i;k_i)=O(1)$ and $D_n=O(1)$ in probability; Further,  the collection of all marginal imbalances  $\big((D_n(i;k_i):w_{m,i}\ne 0, k_i=1,\ldots,m_i,i=1,\ldots, I)\big)_{n\ge 1}$ is  a positive recurrent Markov chain with period $T$.
\item[\rm (c)] When $ w_{m,i}=0$,  then the corresponding marginal imbalance increase  with the rate $\sqrt{n}$, that is, $D_{n,t}(i;k_i)=O(\sqrt{n})$ in probability.
\end{description}
\end{corollary}

\section{Concluding Remarks}
\label{s:conclusion}
\setcounter{equation}{0}
\begin{center}\textbf{}
\end{center}

\vspace{-0.3in}

In this paper we study the theoretical properties of a general family of covariate-adaptive designs. These results provide a unified and fundamental theory about balance properties of covariate-adaptive randomization procedures. In the literature, it is well known that the imbalance is a positive recurrent Markov chain for Efron's (1971) biased coin design (without involving covariates). Markaryan and Rosenberger (2010) studied some exact properties of Efron's biased coin design. Recently, Hu and Hu (2012) showed that the imbalances  $(\bm D_n)_{n\ge 1}$ are positive recurrent Markov chains for a very limited family of covariate-adaptive designs with Efron's bias coin allocation function. The condition (C) in their paper is too restrictive and it is  almost impossible to check this condition in real applications. The results in this paper also provide new insights about  imbalances
on covariate-adaptive randomization procedures: (i) when $w_s>0$ (the within-stratum weight is positive), the imbalances  $(\bm D_n)_{n\ge 1}$ are positive recurrent Markov chains and therefore, all three types of imbalances (within-stratum, marginal and overall) are bounded in probability; (ii) when $w_s=0$ and $w_{m,i}>0$, then the marginal (the $i$-th covariate) and overall imbalances are  bounded in probability, but the within-stratum imbalance is not; (iii) when $w_s=w_{m,i}=0$ for all $i=1,...,I$ and $w_0=1$, only the overall imbalance is bounded in probability.

It is very important to understand statistical inference under covariate-adaptive randomization. In the literature, several authors (Birkett, 1985, Forsythe, 1987,  etc.) have raised concerns about the conservativeness of the unadjusted analysis (such as two-sample t-test) under covariate-adaptive randomization  based on simulation studies. Recently, Shao, Yu and Zhong (2010) studied this problem theoretically under a very special covariate-adaptive biased coin randomization procedure, which is a stratified randomization procedure and only applies to a single covariate case.  Also they focused on a simple homogeneous linear model and only considered a two-sample t-test. This is because the theoretical properties of covariate-adaptive randomization procedures are usually not available in the literature. The results in this paper open a door to study the theoretical behavior of classical statistical inference under covariate-adaptive randomization. For example, based on Corollary 3.2, we can study the behavior of testing hypotheses and other methods under the Pocock and Simon's procedure. We leave these as future research projects.

In this paper, we only consider balancing discrete (categorical) covariates. In the literature, continuous covariates are typically discretized in order to be included in the randomization scheme (Taves, 2010). We may apply the proposed designs to balance continuous covariates after discretized. However, as discussed in Scott et. al. (2002), the breakdown of a continuous covariate into subcategories means increased effort and loss of information. Ciolino et. al. (2011) also pointed out the lack of publicity for practical methods for continuous covariate balancing and lack of knowledge on the cost of failing to balance continuous covariates.
We may consider balancing continuous covariates under similar framework of this paper. However, it is usually difficult to obtain the corresponding theoretical properties. There is not much study in the literature.

The proposed procedures and their properties may be generalized in several ways.  First, we may apply the same idea to problems of unequal ratios (Hu and Rosenberger, 2006). Sometimes, if one treatment is superior (or less costly) than the other, then assigning more patients to the treatment would be more ethical (economical). Second, we may combine the idea in this paper with the ERADE of Hu, Zhang and He (2009) to get a new family of ``covariate-adjusted response-adaptive randomization'' (CARA) procedure (Zhang, Hu, Cheung and Chan, 2007); it could be a real challenge to study the corresponding theoretical properties.  We leave all these as future research topics.

\section{Appendix: Proofs}
\label{s:proofs}
\setcounter{equation}{0}
\begin{center}
\end{center}

\vspace{-0.3in}
Our proofs   are based on the properties of   Markov chains on a countable state space. For general notations and theory for Markov chains we refer to Meyn and Tweedie (1993). For simplification, we write $\bm k= (k_1,\ldots, k_I)$. Let $\Delta\mathscr{D}$ be the state space of $\Delta \bm D_n=\bm D_n-\bm D_{n-1}$, i.e., each $\bm d\in \Delta\mathscr{D}$ has only one non-zero element which is $1$ or $-1$. Also let $\mathscr{F}_n$ be the history $\sigma$-field generated by  the covariates $Z_1,\ldots, Z_n$ and   results of allocation $T_1,\ldots, T_n$.

\smallskip

\noindent{\bf Proof of Proposition \ref{proposition1}.} For (i), taking the summation of $\Lambda_n(\bm k)$ over all $\bm k$ yields
$$ \sum_{\bm k} \Lambda_n(\bm k)=\big(w_s + \sum_{i=1}^I w_{m,i}\prod_{j\ne i} m_j  + w_o m\big) D_n. $$
So $D_n$ is a linear transform of $\bm \Lambda_n$.
Taking the summation of $\Lambda_n(\bm k)$ over all $k_1,\ldots, k_I$ except $k_i$  yields
\begin{align*}
 &\sum_{k_1,\ldots,k_{i-1}, k_{i+1},\ldots, k_I} \Lambda_n(k_1,\ldots,k_I)\\
 =&\big(w_s + w_{m,i}\prod_{j\ne i} m_j)D_n(i;k_i)+\big( \sum_{l\ne i} w_{m,l}\prod_{j\ne i, l} m_j  + w_o \prod_{j\ne i} m_j\big) D_n.
 \end{align*}
Hence, when $w_s+w_{m,i}>0$, each $D_n(i;k_i)$ is a linear transform of $\bm \Lambda_n $ and $D_n$, and so it is  a linear transform of $\bm \Lambda_n$. Finally, when $w_s>0$, it is obvious that each $D_n(\bm k)$ is a linear transform of $\Lambda_n(\bm k)$, $D_n(1;k_1),\ldots, D_n(I;k_I)$ and $D_n$, and so it is a linear transform of $\bm \Lambda_n$. Hence,  when $w_s>0$, $\bm \Lambda_n=\bm L(\bm D_n)$ is a one to one linear map.

\smallskip
For (ii) and (iii), it is sufficient to show the Markov property.  Notice
$$ D_n(\bm k)=D_{n-1}(\bm k)+2\big(T_n-\frac{1}{2}\big)\mathbb{I}\{Z_n=\bm k\}. $$
 Then
\begin{align*}
 \pr(\Delta D_n(\bm k)=1|\mathscr{F}_{n-1})= & g(4\Lambda_{n-1}(\bm k))p(\bm k), \\
 \pr(\Delta D_n(\bm k)=-1|\mathscr{F}_{n-1})= & \big[1-g\big(4\Lambda_{n-1}(\bm k)\big)\big]p(\bm k) \\
 =&g\big(-4\Lambda_{n-1}(\bm k)\big)p(\bm k),\\
\pr(\Delta D_n(\bm k)=0|\mathscr{F}_{n-1})= &1-p(\bm k).
\end{align*}
 For two vectors $\bm x$ and $\bm y$ on $\mathbb{Z}^{m_1\times \ldots\times m_I}$, we write $\bm x\cdot \bm y=\sum_{\bm k}x(\bm k)y(\bm k)$. The conditional probability above can be write in the following form,
 \begin{align}\label{eqprobabDD}
  \pr(\Delta \bm D_n=\bm d|\mathscr{F}_{n-1})=&g(4\bm d\cdot \bm \Lambda_{n-1})|\bm d\cdot \bm p|\\
  =& g\big(4\bm d\cdot \bm L(\bm D_{n-1})\big)|\bm d\cdot \bm p|, \;\; \bm d\in \Delta \mathscr{D},\nonumber
  \end{align}
 which depends only on $\bm \Lambda_{n-1}=\bm L(\bm D_{n-1})$ and is positive. So, conditional on $\bm{D}_{n-1}$, $\bm D_n$ is  conditionally independent of $(\bm{D}_{1},\ldots,\bm{D}_{n-2})$. It follows that $(\bm D_n)_{n\ge 1}$ is a Markov chain on $\mathbb{Z}^{m}$ and is irreducible.

  It is easily seen that
$$ \pr(\Delta (-\bm D_n) =\bm d|\mathscr{F}_{n-1})=g\big(4\bm d\cdot (-\bm \Lambda_{n-1})\big)|\bm d\cdot \bm p|. $$
So, $(-\bm D_n)_{n\ge 1}$ and $(\bm D_n)_{n\ge 1}$  have the same transition probabilities.

For (iii), we consider a more general case. Let $\widetilde{\bm D}=\bm F(\bm D)$ be a linear transform of $\bm D$. We consider the chain $\bm E_n=(\widetilde{\bm D}_n, \bm \Lambda_n)$.  For any point $\bm e$ in the state space $\{(\bm F(\bm d),\bm L(\bm d)):\bm d\in \Delta\mathscr{D}\}$ of $\Delta\bm E_n$,
$$
  \pr(\Delta \bm E_n=\bm e|\mathscr{F}_{n-1})=\sum_{\bm d: (\bm F(\bm d),\bm L(\bm d))=\bm e,\bm d\in \Delta\mathscr{D}} g(4\bm d\cdot \bm \Lambda_{n-1})|\bm d\cdot \bm p|,$$
which depends only on $\bm \Lambda_{n-1}$ and is positive. So, given $\bm E_{n-1}$, $\bm E_n$ is  conditionally independent of $(\bm E_{1},\ldots,\bm E_{n-2})$. It follows that
\begin{equation}\label{eqMarkovProperty}
\Big(\bm E_n=\big(\bm F(\bm D_n), \bm \Lambda_n\big)\Big)_{n\ge 1}\;\;  \text{ is an irreducible Markov chain}.
\end{equation}
 And also, it is easily seen that $(-\bm E_n)_{n\ge 1}$ and $(\bm E_n)_{n\ge 1}$  have the same transition probabilities
 because
  \begin{align*}
  \pr(\Delta \bm E_n=-\bm e|\mathscr{F}_{n-1})=& \sum_{\bm d: (\bm F(\bm d),\bm L(\bm d))=\bm e,\bm d\in \Delta\mathscr{D}} g(-4\bm d\cdot \bm \Lambda_{n-1})|\bm d\cdot \bm p|\\
  =& \sum_{\bm d: (\bm F(\bm d),\bm L(\bm d))=\bm e,\bm d\in \Delta\mathscr{D}} g(4\bm d\cdot (-\bm \Lambda_{n-1}))|\bm d\cdot \bm p|.
  \end{align*}
  The proof of Proposition \ref{proposition1} is completed. $\Box$

\bigskip

\noindent {\bf Proof of Theorem \ref{theorem1}.}
 Define
 $$ V_n=\sum_{\bm k}w_s D_n^2(\bm k)+\sum_{i=1}^I\sum_{k_i=1}^{m_i} w_{m,i} D_n^2(i;k_i)+w_o D_n^2, $$
 $$
\bm D=\left[D(\bm k)\right]_{1\leq k_{1}\leq m_{1},\ldots, 1\leq k_{I}\leq m_{I}},
$$
 and define $\bm\Lambda$ and $V$ in the same way as defining $\bm\Lambda_n$ and $V_n$ with $\bm D$ taking the place of $\bm D_n$.
By Proposition \ref{proposition1} (i), $V_n$ is a function of $\bm \Lambda_n$. Write
$$ V_n=V(\bm\Lambda_n). $$
We will prove the theorem via two steps. First, we will show
that there is bounded set $\mathscr{C}$ and a constant $b$ for which
\begin{equation}\label{eqdriftcondition2}
P_{\lambda} V(\bm\Lambda)- V(\bm \Lambda)\le -1 +b \mathbb{I}_{\bm \Lambda\in \mathscr{C}},
 \end{equation}
  where $P_{\lambda}$ is the transition probability matrix of $\bm \Lambda$,
   $$ P_{\lambda}V(\bm \Lambda)=\sum_{\bm\Lambda^{\prime}\in \bm L(\mathbb{Z}^m)} P_{\lambda}(\bm\Lambda,\bm\Lambda^{\prime}) V(\bm\Lambda^{\prime}), $$
 and $P_{\lambda}(\bm\Lambda,\bm\Lambda^{\prime})$ is the transition probability from state $\bm\Lambda$ to state $\bm\Lambda^{\prime}$.
 In the second step, we will show that  for any integer $r\ge 2$,   there is bounded set $\mathscr{C}$ and a constant $b$ for which
 \begin{equation}\label{eqdriftcondition3}
  P_{\lambda} V^{r+1}(\bm\Lambda)- V^{r+1}(\bm\Lambda)\le - [V(\bm\Lambda)+1]^r   +b \mathbb{I}_{\bm\Lambda\in \mathscr{C}}.
 \end{equation}

 The drift condition (\ref{eqdriftcondition2}) is utilized to show the convergence in probability, and the refined drift condition (\ref{eqdriftcondition3}) is utilized to show the convergence of moments.  In fact, (\ref{eqdriftcondition2}) implies that $(\bm \Lambda_n)_{n\ge 1}$ is a positive (Harris) recurrent Markov chain(c.f., Theorem 11.3.4 of Meyn and Tweedie, 1993), and  so it is bounded in probability and has an invariant probability measure $\pi_{\lambda}$. On the other hand, by (\ref{eqdriftcondition3}) and   Theorem 14.3.7   of Meyn and Tweedie (1993) we conclude that
 $  \pi_{\lambda}[V(\bm\Lambda)+1]^r \le b, $
which implies that
\begin{equation}\label{eqmomentcondition} \sup_n\ep [V(\bm\Lambda_n)+1]^r <\infty
\end{equation}
by
Theorem 14.3.6 of Meyn and Tweedie (1993).  Notice that
 \begin{align*}
& |\Lambda_{n}(\bm k)|^2\le \big(w_{o}|D_{n}|+\sum_{i=1}^{I}w_{m,i}|D_{n}(i;k_{i})|+w_{s}|D_n(\bm k)|\big)^2\\
& \le  \big(w_{o}|D_{n}|^2+\sum_{i=1}^{I}w_{m,i}|D_{n}(i;k_{i})|^2+w_{s}|D_{n}(\bm k)|^2\big)(w_{o} +\sum_{i=1}^{I}w_{m,i} +w_{s} ),\\
&=  w_{o}|D_{n}|^2+\sum_{i=1}^{I}w_{m,i}|D_{n}(i;k_{i})|^2+w_{s}|D_{n}(\bm k)|^2,
\end{align*}
which implies that $ \|\bm \Lambda_n\|^2\le m V(\bm\Lambda_n)$. It follows that
$ \sup_n\ep\|\bm \Lambda_n\|^{2r}<\infty. $
Thus, we conclude that $(\bm \Lambda_n)_{n\ge 1}$ is a positive recurrent Markov chain with $\ep\|\bm\Lambda_n\|^r=O(1)$ for all $r>0$. (i)-(iii) follows from   Proposition \ref{proposition1} (i).

\bigskip
Now, we begin the proofs of (\ref{eqdriftcondition2}) and (\ref{eqdriftcondition3}).   Given $Z_n=\bm k$, if $T_n=1$, then
 \begin{align*}
 V_n-V_{n-1}=& w_s \left\{ \big(D_{n-1}(\bm k)+1\big)^2-D_{n-1}^2(\bm k)\right\}\\
 &+\sum_{i=1}^I w_{m,i} \left\{ \big(D_{n-1}(i;k_i)+1\big)^2-D_{n-1}^2(i;k_i)\right\}\\
 &+  w_o \left\{ (D_{n-1}+1)^2-D_{n-1}^2 \right\}\\
 =& 2\Lambda_{n-1}(\bm k)+1,
 \end{align*}
 while, if $T_n=0$, then
 $V_n-V_{n-1}= -2\Lambda_{n-1}(\bm k)+1. $
 So,
\begin{equation}\label{eqV} V_n-V_{n-1}=4\Lambda_{n-1}(\bm k)\left(T_n-\frac{1}{2}\right)\mathbb{I}\{Z_n=\bm k\}+1.
\end{equation}
 Hence
\begin{align*}
  \ep\left[V_n-V_{n-1}\big|Z_n=\bm k,\mathscr{F}_{n-1}\right]
 =&2\Lambda_{n-1}(\bm k)\left[ g(4\Lambda_{n-1}(\bm k))-g(-4\Lambda_{n-1}(\bm k))\right]+1\\
 =& -4\Big|\Lambda_{n-1}(\bm k)\Big|\cdot\left( \frac{1}{2}-g\left(4|\Lambda_{n-1}(\bm k)|\right)\right)+1.
 \end{align*}
 It follows that
 \begin{align}\label{eqcondE}
  \ep\left[V_n\big|\mathscr{F}_{n-1}\right]-V_{n-1}
  = -4S(\bm\Lambda_{n-1})+1,
 \end{align}
 where $$ S(\bm\Lambda_{n-1})=\sum_{\bm k}\Big|\Lambda_{n-1}(\bm k)\Big|\cdot\left( \frac{1}{2}-g\left(4|\Lambda_{n-1}(\bm k)|\right)\right)p(\bm k)$$
 is a nonnegative function of $\bm \Lambda_{n-1}$ since $g(x)\le 1/2$ when $x\ge 0$.

 Recall that $\bm D_n$ and $\bm\Lambda_n$ are   irreducible Markov chains with period two on $\mathbb{Z}^m$ and $\bm L(\mathbb{Z}^m)$, respectively.     Notice that $V_n=V(\bm\Lambda_n)$ is a non-negative function of $\bm\Lambda_n$. Equation (\ref{eqcondE}) tells us that the drift function of the Markov chain $\bm\Lambda_n$ is
 \begin{align*}
  \ep\left[V(\bm\Lambda_n)\big|\bm \Lambda_{n-1}\right]-V(\bm\Lambda_{n-1})
 =\ep[V_n-V_{n-1}|\bm \Lambda_1,\ldots,\bm \Lambda_{n-1}]= 1-4 S(\bm \Lambda_{n-1})
 \end{align*}
 by the Markov-property. That is
 $$P_{\lambda}V(\bm \Lambda)-V(\bm\Lambda)=1-4 S(\bm \Lambda).$$
  Next, we need to check the drift-criteria condition (\ref{eqdriftcondition2}).   It is sufficient to show that
 \begin{equation}\label{eq3}
  \bm \Lambda\text{ is bounded} \Leftrightarrow S(\bm\Lambda) \text{ is bounded}.
  \end{equation}

 Write $\underline{p}=\min_{\bm k}p(\bm k)$. By the assumption that $\limsup_{x\to +\infty} g(x)<1/2$, there exist two positive constants $\underline{\lambda}$ and $\underline{g}$ such that
 $ 1/2- g(4x)\ge \underline{g}>0$  for all  $x\ge \underline{\lambda}$.
  Then
  \begin{align*}
   S(\bm\Lambda)=&\sum_{\bm k}\big|\Lambda(\bm k)\big|\cdot\left( \frac{1}{2}-g\left(4|\Lambda(\bm k)|\right)\right)p(\bm k)
   \ge    \underline{g}\cdot \sum_{\bk k:|\Lambda(\bm k)|\geq \underline{\lambda}}\big|\Lambda(\bm k)\big|p(\bm k)\\
   =&\underline{g}\cdot\Big( \sum_{\bk k}\big|\Lambda(\bm k)\big|p(\bm k)-  \sum_{\bk k:|\Lambda(\bm k)|< \underline{\lambda}}\big|\Lambda(\bm k)\big|p(\bm k)\Big) \\
   \ge &\underline{g}\cdot\Big( \sum_{\bk k}\big|\Lambda(\bm k)\big|p(\bm k)- \underline{\lambda}\cdot\sum_{\bk k:|\Lambda(\bm k)|< \underline{\lambda}}  p(\bm k)\Big)
   \ge  \underline{g}\cdot\Big( \sum_{\bk k}\big|\Lambda(\bm k)\big|p(\bm k)- \underline{\lambda}\Big).
  \end{align*}
   It follows that
  $$ \underline{g}\cdot \underline{p}\cdot\sum_{\bm k}\big|\Lambda(\bm k)\big|-\underline{\lambda}\le S(\bm\Lambda)\le \sum_{\bm k }\big|\Lambda(\bm k)\big|. $$
  (\ref{eq3}) is proved.  From  (\ref{eq3}), it follows that there is bounded set $\mathscr{C}$ and a constant $b$ such that  the drift condition (\ref{eqdriftcondition2}) is satisfied.

For  verifying (\ref{eqdriftcondition3}), we shall refine drift condition (\ref{eqdriftcondition2}).   For given $Z_n=\bm k$, by (\ref{eqV}) we have
\begin{align*}
 V_n=&V_{n-1}+1+4\Lambda_{n-1}(\bm k)(T_n-\frac{1}{2})\mathbb{I}\{Z_n=\bm k\}\\
 \widehat{=}& V_{n-1}+1+\xi .
 \end{align*}
It is obvious that
$$|\xi|=2|\Lambda_{n-1}(\bm k) |\le 2 \sqrt{ V_{n-1}},\;\; \ep[\xi|\mathscr{F}_{n-1}]=-4S(\bm\Lambda_{n-1}). $$
It follows that
\begin{align*}
V_n^{r+1}- & V_{n-1}^{r+1}=(r+1)(V_{n-1}+1)^r\xi\\
&\quad +\left\{(V_{n-1}+1)^{r+1}-V_{n-1}^{r+1}+\sum_{k=2}^{r+1}\binom{r+1}{k}\xi^k(V_{n-1}+1)^{r+1-k}\right\}\\
&\qquad  \le (r+1)(V_{n-1}+1)^r\xi+C_r (V_{n-1}+1)^r,
\end{align*}
where $C_r$ is a constant which depends on $r$. It follows that
\begin{align*} \ep[V_n^{r+1}|\mathscr{F}_{n-1}]-V_{n-1}^{r+1}\le -4 (r+1)(V_{n-1}+1)^rS(\bm\Lambda_{n-1})+C_r (V_{n-1}+1)^r.
\end{align*}
That is
$$ P_{\lambda} V^{r+1}(\bm\Lambda)- V^{r+1}(\bm\Lambda)\le [V(\bm\Lambda)+1]^r\left\{-4(r+1)S(\bm\Lambda)+C_r\right\}, $$
which, together with (\ref{eq3}), implies  (\ref{eqdriftcondition3}). The proof of Theorem \ref{theorem1} is now completed. $\Box$

\bigskip
\noindent{\bf Proof of Theorem \ref{theorem2}.}
We first prove that
\begin{equation}\label{eqmomentbound}
\sup_n\ep\left|\frac{ D_n(\bm k)}{\sqrt{n}}\right|^r<\infty, \;\; \forall \bm k, \; r>0.
  \end{equation}
 Notice
$$ D_n(\bm k)= D_{n-1}(\bm k)+2\Big(T_n-\frac{1}{2}\Big)\mathbb{I}\{\bm Z_n=\bm k\}. $$
For simplification, we write   $\overline{g}_{n-1,\bm k}=g\big(4\Lambda_{n-1}(\bm k)\big)-g\big(-4\Lambda_{n-1}(\bm k)\big)$ and
$\overline{g}_{\bm k}=g\big(4\Lambda(\bm k)\big)-g\big(-4\Lambda(\bm k)\big)$. It follows that
\begin{align*} \ep[D_n(\bm k)|\mathscr{F}_{n-1}]= & D_{n-1}(\bm k) +\left[g\big(4\Lambda_{n-1}(\bm k)\big)-g\big(-4\Lambda_{n-1}(\bm k)\big)\right]p(\bm k)\\
=& D_{n-1}(\bm k) + \overline{g}_{n-1,\bm k}p(\bm k).
\end{align*}
So
\begin{equation}\label{eqDnij}D_n(\bm k)=\sum_{l=1}^n \big(D_l(\bm k)-\ep[D_l(\bm k)|\mathscr{F}_{l-1}]\big)+
 p(\bm k)\sum_{l=0}^{n-1}\overline{g}_{l,\bm k}.
\end{equation}
The first term on the right hand above is $O(\sqrt{n})$ in $L_r$,
because $\{D_n(\bm k)-\ep[D_n(\bm k)|\mathscr{F}_{n-1}]\}$ is a sequence of bounded martingale differences.
Now, we consider the second term. Notice $\overline{g}_{\bm k}$ is bounded by $1$.
By (\ref{eqdriftcondition2}) and Theorem 17.4.2 of Meyn and Tweedie (1993), there is a constant $R$ such that the Poisson's equation
\begin{equation}\label{PossionEquation}
 \widehat{g}-P_{\lambda}\widehat{g}=\overline{g}_{\bm k}-\pi_{\lambda}\overline{g}_{\bm k}
 \end{equation}
has a solution $\widehat{g}=\widehat{g}_{\bm k}=\widehat{g}_{\bm k}(\bm\Lambda)$ which is a function of $\bm\Lambda$ defined on the state space of $\bm \Lambda$  with $|\widehat{g}|\le R(V+1)$, where $P_{\lambda}$ is the transition probability matrix of $\bm\Lambda$.
On the other hand,  notice that $(-\bm \Lambda_n)_{n\ge 1}$ and $(\bm \Lambda_n)_{n\ge 1}$ have the same transition probabilities.
It follows that the invariant probability measure $\pi_{\lambda}$ is symmetric, that is, under $\pi_{\lambda}$,  $(-\bm \Lambda_n)_{n\ge 1}$ and $(\bm \Lambda_n)_{n\ge 1}$  are identically distributed. It follows that $\pi_{\lambda}\overline{g}_{\bm k}=\pi_{\lambda}[g(4\Lambda(\bm k))-g(-4\Lambda(\bm k))]=0$. So,
\begin{equation}\label{PossionEquation2}
 \widehat{g}_{\bm k}-P_{\lambda}\widehat{g}_{\bm k}=\overline{g}_{\bm k}.
 \end{equation}
Now, write $\widehat{g}_n =\widehat{g}_{\bm k}(\bm\Lambda_n)$. It follows that
\begin{align*}
\sum_{l=0}^{n-1}\overline{g}_{l, \bm k}=& \sum_{l=0}^{n-1}\left\{\widehat{g}_l-P_{\lambda}\widehat{g}_l\right\}= \sum_{l=0}^{n-1}\left\{\widehat{g}_l-\ep[\widehat{g}_{l+1}|\mathscr{F}_l]\right\}\\
=& \sum_{l=0}^{n-1}\left\{\widehat{g}_l-\ep[\widehat{g}_l|\mathscr{F}_{l-1}]\right\}
+\ep[\widehat{g}_1]-\ep[\widehat{g}_n|\mathscr{F}_{n-1}].
\end{align*}
Hence for any $r\ge 1$,
\begin{align*}
&\ep\left|\frac{1}{\sqrt{n}}\sum_{l=0}^{n-1} \overline{g}_{l,\bm k}\right|^{2r}\\
\le &C \ep\left|\frac{1}{n}\sum_{l=0}^{n-1} \ep\left[\left\{\widehat{g}_l-\ep[\widehat{g}_l|\mathscr{F}_{l-1}]\right\}^2\big|\mathscr{F}_{l-1}\right]\right|^{r}+C \frac{\ep\widehat{g}_1^{2r}+\ep\widehat{g}_{n}^{2r}}{n^r} \\
\le & C\frac{\sum_{l=0}^{n}\ep\widehat{g}_l^{2r}}{n}\le C\sup_l\ep  \widehat{g}_l^{2r}\le CR^{2r} \sup_n\ep(V_n+1)^{2r}<\infty
\end{align*}
 by (\ref{eqmomentcondition}). (\ref{eqmomentbound}) is now proved.

\smallskip
Next we  prove (iv). We prove (\ref{varainceofD1}) first.   Fix $\bm k=(k_1,\ldots,k_I)$. Let $\widehat{g}=\widehat{g}_{\bm k}$ be the solution of the Poisson's equation (\ref{PossionEquation}) which is a function of $\bm\Lambda$. Let $\bm B_{\bm k}\in \bm L(\Delta\mathscr{D})$ be the element whose value is $\bm \Lambda_n-\bm\Lambda_{n-1}$ with $D_n(\bm k)-D_{n-1}(\bm k)=1$. That is, the $\bk l$-th element of $\bm B_{\bm k}$ is
  $w_s\mathbb{I}\{\bm l=\bm k\}+\sum_{i=1}^I w_{m,i}\mathbb{I}\{l_i=k_i\}+w_o$. We will show that (\ref{varainceofD1}) holds with
   \begin{equation}\label{varainceofD2}
  \sigma^2(\bm k)=\pi_{\lambda}[h_{\bm k,\bm k}(\bm \Lambda)],
 \end{equation}
 where
 $$ h_{\bm k,\bm k}(\bm\Lambda)=p(\bm k)+2p^2(\bm k)\left[\widehat{g}_{\bm k}(\bm\Lambda+\bm B_{\bm k})g(4\Lambda(\bm k))-\widehat{g}_{\bm k}(\bm\Lambda-\bm B_{\bm k})g(-4\Lambda(\bm k))\right]. $$
 Denote
 \begin{equation}
 \label{eqMdifference}
 \Delta M_{n,\bm k}=D_n(\bm k)-D_{n-1}(\bm k)+p(\bm k)[\widehat{g}_{\bm k}( \bm \Lambda_n)-\widehat{g}_{\bm k}(\bm \Lambda_{n-1})].
 \end{equation}
  Then
 \begin{align*}
 \ep[\Delta M_{n,\bm k}|\mathscr{F}_{n-1}]=&\ep[D_n(\bm k)|\mathscr{F}_{n-1}]- D_{n-1}(\bm k) +p(\bm k)[P_{\lambda}\widehat{g}_{\bm k}( \bm \Lambda_{n-1})-\widehat{g}_{\bm k}(\bm \Lambda_{n-1})]\\
 =&\ep[D_n(\bm k)|\mathscr{F}_{n-1}]- D_{n-1}(\bm k)-p(\bm k)\overline{g}_{n-1,\bm k}=0.
 \end{align*}
 So, $\{\Delta M_{n,\bm k}\}$ is a sequence of martingale differences with
\begin{equation}\label{martingale}
M_{n,\bm k}=\sum_{l=1}^n \Delta M_{l,\bm k}=D_n(\bm k)-D_0(\bm k)+p(\bm k)[\widehat{g}_{\bm k}( \bm \Lambda_n)-\widehat{g}_{\bm k}(\bm \Lambda_0)]
 \end{equation}
and
\begin{align*}
&\ep[(\Delta M_{n,\bm k})^2|\mathscr{F}_{n-1}]\\
=&\ep[ (\Delta D_n(\bm k))^2|\mathscr{F}_{n-1}]+2p(\bm k) \ep\left[[\widehat{g}_{\bm k}(\bm \Lambda_n)-\widehat{g}_{\bm k}(\bm \Lambda_{n-1})](\Delta D_n(\bm k))\big|\mathscr{F}_{n-1}\right] \\
&+p^2(\bm k)\ep \left[[\widehat{g}_{\bm k}( \bm \Lambda_n)-\widehat{g}_{\bm k}(\bm \Lambda_{n-1})]^2\big|\mathscr{F}_{n-1}]\right]\\
=&\ep[ (\Delta D_n(\bm k))^2|\mathscr{F}_{n-1}]+2p(\bm k) \ep[\widehat{g}_{\bm k}(\bm \Lambda_n)(\Delta D_n(\bm k))|\mathscr{F}_{n-1}] \\
&-2p(\bm k)\widehat{g}_{\bm k}(\bm \Lambda_{n-1}) \ep[\Delta D_n(\bm k)|\mathscr{F}_{n-1}]\\
&-2p^2(\bm k)\widehat{g}_{\bm k}(\bm \Lambda_{n-1}) \left[\ep[\widehat{g}_{\bm k}(\bm \Lambda_n)|\mathscr{F}_{n-1}]-\widehat{g}_{\bm k}(\bm \Lambda_{n-1})\right]\\\
&+p^2(\bm k) \left(\ep[\widehat{g}_{\bm k}^2(\bm\Lambda_n)|\mathscr{F}_{n-1}]-\widehat{g}_{\bm k}^2(\bm\Lambda_{n-1})\right)\\
=&h_{\bm k,\bm k}(\bm\Lambda_{n-1}) -2p^2(\bm k)\widehat{g}_{\bm k}(\bm \Lambda_{n-1}) \overline{g}_{n-1}(\bm k) \\
&-2p^2(\bm k)\widehat{g}_{\bm k}(\bm \Lambda_{n-1}) \left(P_{\lambda}\widehat{g}_{\bm k}(\bm \Lambda_{n-1}) -\widehat{g}_{\bm k}(\bm \Lambda_{n-1})\right)\\
&+p^2(\bm k) \left(\ep[\widehat{g}_{\bm k}^2(\bm\Lambda_n)|\mathscr{F}_{n-1}]-\widehat{g}_{\bm k}^2(\bm\Lambda_{n-1})\right)\\
=&h_{\bm k,\bm k}(\bm\Lambda_{n-1})  +p^2(\bm k) \left(\ep[\widehat{g}_{\bm k}^2(\bm\Lambda_n)|\mathscr{F}_{n-1}]-\widehat{g}_{\bm k}^2(\bm\Lambda_{n-1})\right),
\end{align*}
where the last equality is  due to the equation (\ref{PossionEquation2}). It follows that
$$ \ep M_{n,\bm k}^2=\sum_{l=0}^{n-1} \ep h_{\bm k,\bm k}(\bm\Lambda_l)+p^2(\bm k) \left(\ep[\widehat{g}_{\bm k}^2(\bm\Lambda_n)]-\ep[\widehat{g}_{\bm k}^2(\bm\Lambda_0)]\right).$$
For $h_{\bm k,\bm k}(\cdot)$, it is easily seen that
$ \pi_{\lambda}[h_{\bm k,\bm k}(\bm\Lambda)]=\sigma^2(\bm k)\ge 0 $
and
\begin{align*}
 |h_{\bm k,\bm k}(\bm\Lambda)|\le & 1+2|\widehat{g}_{\bm k}(\bm\Lambda\pm \bm B_{\bm k})|\\
 \le & 1+2R(V(\bm\Lambda\pm \bm B_{\bm k})+1)\le c_0(V(\bm\Lambda)+1).
 \end{align*}
By (\ref{eqdriftcondition3}) (where $r=1$) and applying Theorem 17.4.2 of Meyn and Tweedie (1993) again, we have a function $\widehat{h}(\bm\Lambda)$ such that
$$ \widehat{h}-P_{\lambda}\widehat{h}=h_{\bm k,\bm k}- \pi_{\lambda}[h_{\bm k,\bm k}] \;\;\text{ and } \;\; |\widehat{h}|\le c (V^2+1). $$
It follows that
\begin{align}\label{eqVaRofM}
&\ep[M_{n,\bm k}^2]=\sum_{l=0}^{n-1}\pi_{\lambda}[h_{\bm k,\bm k}(\bm \Lambda_l)]+\sum_{l=0}^{n-1} \ep\left\{\widehat{h}(\bm\Lambda_l)-P_{\lambda}\widehat{h}(\bm\Lambda_l)\right\} \nonumber\\
&\qquad \qquad \quad +p^2(\bm k) \left(\ep[\widehat{g}_{\bm k}^2(\bm\Lambda_n)]-\ep[\widehat{g}_{\bm k}^2(\bm\Lambda_0)]\right)\nonumber\\
& =n\sigma^2(\bm k)+\sum_{l=0}^{n-1} \left\{\ep\widehat{h}(\bm\Lambda_l)-\ep\widehat{h}(\bm\Lambda_{l+1})\right\}+p^2(\bm k) \left(\ep[\widehat{g}_{\bm k}^2(\bm\Lambda_n)]-\ep[\widehat{g}_{\bm k}^2(\bm\Lambda_0)]\right)\nonumber\\
& =n\sigma^2(\bm k)+  \left\{\ep\widehat{h}(\bm\Lambda_0)-\ep\widehat{h}(\bm\Lambda_n)\right\}+p^2(\bm k) \left(\ep[\widehat{g}_{\bm k}^2(\bm\Lambda_n)]-\ep[\widehat{g}_{\bm k}^2(\bm\Lambda_0)]\right)\nonumber\\
& =n\sigma^2(\bm k)+O(1),
\end{align}
by notice that $\ep V^2(\bm\Lambda_n)$ is bounded and $|\widehat{g}_{\bm k}|\le R(V+1)$. Hence (\ref{varainceofD1}) is proved.

Notice $\ep|\Delta M_{n,\bm k}|^r\le c+c\sup_n\ep|\widehat{g}_{\bm k}(\bm\Lambda_n)|^r\le c+c \sup_n\ep V^r(\bm\Lambda_n)<\infty$
by (\ref{eqmomentcondition}). By the central limit theorem for martingales, we conclude that
$$ \frac{D_n(\bm k)}{\sqrt{n}}=\frac{M_{n,\bm k}+O_P(1)}{\sqrt{n}}\overset{d}\to N(0,\sigma^2(\bm k)). $$
The asymptotic normality  (\ref{eqCLTofD}) is proved. The asymptotic normality, together with (\ref{eqmomentbound}), implies (\ref{eqth1.3}).  The proof of (iv) is completed.

\smallskip

For (v) and (vi), notice that if $D_n(\bm k)=o_P(\sqrt{n})$ then  $D_n(\bm k)=o(\sqrt{n})$ in $L_2$ by (\ref{eqmomentbound}). So, $\sigma(\bm k)=0$. Hence $\ep D_n^2(\bm k)=O(1)$ by (\ref{varainceofD1}).  (v) is proved.

Further, $\ep M_{n,\bm k}^2=O(1)$ by  (\ref{eqVaRofM}). By the Martingale Convergence Theorem, there is a random variable $M_{\infty}$ such that
\begin{equation}\label{eqconvergM}
M_{n,\bm k} \to M_{\infty}\;\; a.s. \text{ and }\;\; M_{n,\bm k}=\ep[M_{\infty}|\mathscr{F}_n].
\end{equation}
On the other hand,  the sequence $\big( (D_n(\bm k), \bm\Lambda_n)\big)_{n\ge 1}$ is an irreducible Markov chain by  (\ref{eqMarkovProperty}). Hence it is a  positive (Harris) recurrent Markov chain by Proposition 18.3.1 of Meyn and Tweedie (1993) due to the fact that it is bounded in probability. Recall the equation (\ref{martingale}). The left hand is a martingale which is convergent almost surely due to (\ref{eqconvergM}), while, the right hand is a function of a  positive  (Harris) recurrent Markov chain.  It follows that the limit $M_{\infty}$ must be a constant. So
$$M_{n,\bm k}=\ep[M_{\infty}|\mathscr{F}_n]=const\;\; a.s. $$
It is obvious that $M_0=0$. Hence
$M_{n,\bm k}\equiv 0\;\; a.s. $
It follows that
\begin{equation}\label{eqfunctionDk}
D_n(\bm k)=D_0(\bm k)-p(\bm k)[\widehat{g}_{\bm k}( \bm \Lambda_n)-\widehat{g}_{\bm k}(\bm \Lambda_0)],
\end{equation}
which implies that $D(\bm k)$ is a function of $\bm \Lambda$. Up to now, we arrive at the conclusion that if $D_n(\bm k)=o_P(\sqrt{n})$, then $D(\bm k)$ is a function of $\bm \Lambda$.

 Finally, we show that (\ref{eqfunctionDk}) is a contradiction when $w_s=0$. Choose a stratum $\bm k^{\ast}$ such that $k^{\ast}_i\ne k_i$, $i=1,\ldots, I$.   Recall that $\bm B_{\bm j}$ is the value of $\Delta \bm\Lambda_n$ when $\Delta D_n(\bm j)=1$ and the $\bm l$-th value is $\bm B_{\bm j}(\bm l) =\sum_{i=1}^Iw_{m,i}\mathbb{I}\{l_i=j_i\}+w_o$. It follows that
\begin{align*}
\bm B_{\bm k} (\bm l)+\bm B_{\bm k^{\ast}}(\bm l) =&\sum_{i=1}^Iw_{m,i}\mathbb{I}\{l_i=k_i\text{ or } k_i^{\ast}\}+2w_o\\
=& \bm B_{k_1^{\ast},k_2,\ldots,k_I}(\bm l)+\bm B_{k_1,k_2^{\ast},\ldots,k_I^{\ast}}(\bm l).
\end{align*}
That is
$$ \bm B_{\bm k} +\bm B_{\bm k^{\ast}}
-\bm B_{k_1^{\ast},k_2,\ldots,k_I}-\bm B_{k_1,k_2^{\ast},\ldots,k_I^{\ast}} =\bm 0.$$
 It follows that on the event
 $E=\{\Delta D_n(\bm k)=1,\Delta D_{n+1}(\bm k^{\ast})=1,
\Delta D_{n+2}(k_1^{\ast},k_2,\ldots,k_I)=-1, \Delta D_{n+3}(k_1,k_2^{\ast},\ldots,k_I^{\ast})=-1\}$, the value of $\bm \Lambda$ does not change, and so the value of the right hand of (\ref{eqfunctionDk}) does not change, i.e.,
$$-p(\bm k)[ \widehat{g}_{\bm k}( \bm \Lambda_{n+3})-\widehat{g}_{\bm k}(\bm \Lambda_{n-1})]=0. $$
However, on the event $E$, $D_{n+3}(\bm k)-D_{n-1}(\bm k)=1$. On the other hand, it is easily seen that, conditional on $\bm D_{n-1}$, the probability of $E$ is positive.  We get a   contradiction to (\ref{eqfunctionDk}). The proof of Theorem \ref{theorem2} is now completed. $\Box$

\bigskip
{\bf Proof of Theorem \ref{theorem3}.} For the martingale difference in (\ref{eqMdifference}), by (\ref{PossionEquation2}) we can also show that
 \begin{align*}
 &\ep[\Delta M_{n,\bm k}\Delta M_{n,\bm l} |\mathscr{F}_{n-1}]\\
 =& p(\bm k)\ep[\widehat{g}_{n,\bm k}\Delta D_n(\bm l)|\mathscr{F}_{n-1}]+p(\bm l)\ep[\widehat{g}_{n,\bm l}\Delta D_n(\bm k)|\mathscr{F}_{n-1}]\\
 &+ p(\bm k)p(\bm l)\left[ \ep[\widehat{g}_{n,\bm k}\widehat{g}_{n,\bm l}|\mathscr{F}_{n-1}]-\widehat{g}_{n-1,\bm k}\widehat{g}_{n-1,\bm l}\right]\\
 =:& h_{\bm k,\bm l}(\bm \Lambda_{n-1})+p(\bm k)p(\bm l)\left[ \ep[\widehat{g}_{n,\bm k}\widehat{g}_{n,\bm l}|\mathscr{F}_{n-1}]-\widehat{g}_{n-1,\bm k}\widehat{g}_{n-1,\bm l}\right]
 \end{align*}
for $\bm k\ne \bm l$. With the same argument as (\ref{eqVaRofM}), we have
$$ \ep[ M_{n,\bm k}  M_{n,\bm l}  ]=\sum_{l=1}^n\ep[\Delta M_{l,\bm k}\Delta M_{l,\bm l} ]
=n \pi_{\lambda}[h_{\bm k,\bm l}] +O(1). $$
Write $\sum_{\bm k\setminus k_i}$ for taking the summation over all $k_1,\ldots,k_{i-1},k_{i+1},\ldots,k_I$. Define
$ M_{n}(i;k_i)= \sum_{\bm k\setminus k_i}  M_{n,\bm k}. $
It follows that
$$ \ep[M_n^2(i;k_i)]=n \sigma^2(i;k_i)+O(1). $$
Taking the summation on both side of (\ref{martingale}) over $\bm k\setminus k_i$ yields
$$ M_n(i;k_i)=D_n(i;k_i)-D_0(i;k_i)+g_{i;k_i}(\bm\Lambda_n) $$
where $g_{i;k_i}(\bm\Lambda_n)=\sum_{\bm k\setminus k_i}p(\bm k)[\widehat{g}_{\bm k}( \bm \Lambda_n)-\widehat{g}_{\bm k}(\bm \Lambda_0)]$ is a function of $\bm\Lambda_n$. Hence
\begin{equation}\label{eqVaRofMMarg} \ep[D_n^2(i;k_i)]=n \sigma^2(i;k_i)+O(\sqrt{n} \sigma(i;k_i)).
\end{equation}
With the same argument as showing (\ref{eqfunctionDk}),  if  $\sigma^2(i;k_i)=0$, then we will have $M_n(i;k_i)\equiv 0$ and
\begin{equation}\label{eqfunctionDiki}
D_n(i;k_i)-D_0(i;k_i)=-g_{i;k_i}(\bm\Lambda_n).
\end{equation}
Under the condition $w_s+w_{m,i}=0$, $\bm \Lambda_n$ is a linear transform of $\big(D_n(j;k_j): j=1,\ldots,i-1,i+1,\ldots, I, k_j=1,\ldots, m_j\big)$ which excludes the values of marginal imbalances $D(i;l_i),l_i=1,\ldots, m_i$, of the $i$th covariate. It follows that for $\Delta D_n(\bm k)=1$ and
$\Delta D_n(k_1,\ldots,k_{i-1},k_i^{\ast}, k_{i+1},\ldots, k_I)=1$, $\Delta \bm\Lambda_n$ is the same, and so the values changed  in the right hand of (\ref{eqfunctionDiki}) are the same. Obviously, the values changed in the left hand of (\ref{eqfunctionDiki}) are different, one of which is $1$, the other is $0$. We get a contradiction. $\Box$

\bigskip
{\bf Proof of Theorem \ref{theorem1mul}.} The basic idea of this proof is similar to that of Theorems \ref{theorem1}-\ref{theorem3}. Here we only give the difference. At first, it can be verified that, for the stratum $\bm k^{\ast}$,
$ Imb_{n-1,t}-Imb_{n-1,h}=2\big(\Lambda_{n-1,t}(\bm k^{\ast})-\Lambda_{n-1,h}(\bm k^{\ast})\big). $
So the order of $Imb_{n-1,(1)}\le \ldots\le Imb_{n-1,(T)}$ and the order of $\Lambda_{n-1,(1)}(\bm k^{\ast})\le\ldots \Lambda_{n-1,(T)}(\bm k^{\ast})$ are the same  and therefore the allocation probabilities of the $n$-th patient are  functions of $\bm\Lambda_{n-1}$. Define
$$ V_n=\sum_{h=1}^T\Big\{w_s\sum_{\bm k}D^2_{n,h}(\bm k)+\sum_{i=1}^I\sum_{k_i=1}^{m_i}w_{m,i}D_{n,h}^2(i;k_i)
+w_0D_{n,h}^2\Big\}. $$
Then $V_n$ is a norm-like function of $\bm \Lambda_n$ as the form $V_n=V(\bm \Lambda_n)$. Given  $\bm Z_n=\bm k$,  it can be verified that
$$ V_n-V_{n-1}=2\Lambda_{n-1,(t)}(\bm k)+\frac{T-1}{T}$$
if the patient is allocated to treatment $(t)$.
It follows that
$$ \ep[V_n|\mathscr{F}_{n-1}]-V_{n-1}=-2S(\bm \Lambda_{n-1})+\frac{T-1}{T}, $$
where
$  S(\bm \Lambda_{n-1})=-\sum_{\bm k} \sum_{t=1}^T p_t \cdot \Lambda_{n-1,(t)}(\bm k)\cdot p(\bm k). $
Note that $\sum_{t=1}^T  \Lambda_{n-1,(t)}(\bm k)=0$ and $(p_t-p_h)[ \Lambda_{n-1,(t)}(\bm k)- \Lambda_{n-1,(h)}(\bm k)]\le 0$. We have
\begin{align*}
&  2T\sum_{t=1}^T p_t\cdot \Lambda_{n-1,(t)}(\bm k)\\
= & 2T  \sum_{t=1}^T p_t \cdot \Lambda_{n-1,(t)}(\bm k)-2\big[  \sum_{t=1}^T p_t\big]\big[\sum_{t=1}^T \Lambda_{n-1,(t)}(\bm k)\big]\\
 = &\sum_{t,h=1}^T(p_t-p_h)[ \Lambda_{n-1,(t)}(\bm k)- \Lambda_{n-1,(h)}(\bm k)] \\
 \le & -(p_1-p_T)[ \Lambda_{n-1,(T)}(\bm k)- \Lambda_{n-1,(1)}(\bm k)]\\
\le & -(p_1-p_T)\frac{1}{T}\sum_t\big|\Lambda_{n-1,t}(\bm k)\big|.
\end{align*}
It follows that $S(\bm\Lambda_{n-1})\ge 0$ and
\begin{align*}
 \frac{(p_1-p_T)\min_{\bm k}p(\bm k)}{2T^2}\sum_{\bm k} \sum_{t=1}^T \big|\Lambda_{n-1,t}(\bm k)\big|
\le   S(\bm \Lambda_{n-1})\le \sum_{\bm k} \sum_{t=1}^T \big|\Lambda_{n-1,t}(\bm k)\big|,
\end{align*}
which implies that
$$ \bm \Lambda\text{ is bounded} \Leftrightarrow S(\bm\Lambda) \text{ is bounded}. $$
Further, if we denote
$$ V_n-V_{n-1}=\xi+\frac{T-1}{T},$$
Then
$$ |\xi|\le 2\sqrt{V_{n-1}}, \; \ep [\xi|\mathscr{F}_{n-1}]=-2S(\bm\Lambda_{n-1}),  $$
and for any integer $r\ge 2$,
\begin{align*}
V_n^{r+1}- & V_{n-1}^{r+1}= \Big(V_{n-1}+\frac{T-1}{T}+\xi\Big)^{r+1}\\
&\qquad  \le (r+1)\Big(V_{n-1}+\frac{T-1}{T}\Big)^r\xi+C_r \Big(V_{n-1}+\frac{T-1}{T}\Big)^r.
\end{align*}
  It follows that
\begin{align*} \ep[V_n^{r+1}|\mathscr{F}_{n-1}]-V_{n-1}^{r+1}\le  & -2 (r+1)\Big(V_{n-1}+\frac{T-1}{T}\Big)^rS(\bm\Lambda_{n-1})\\
& +C_r \Big(V_{n-1}+\frac{T-1}{T}\Big)^r
\end{align*}
We arrive at the drift condition:
$$ P_{\lambda} V^{r+1}(\bm\Lambda)- V^{r+1}(\bm\Lambda)\le \Big[V(\bm\Lambda)+\frac{T-1}{T}\Big]^r\Big\{-2(r+1)S(\bm\Lambda)+C_r\Big\}. $$
It follows that $(\bm \Lambda_n)_{n\ge 1}$ is a positive recurrent Markov chain with $\sup_n \ep V^r(\bm \Lambda_n)<\infty$.
The first part of the conclusions in the theorem is proved and (i)-(iii) follow.

For (iv) and (v), we let $\pi_{\lambda}$ be the  invariant probability measure of the positive recurrent Markov chain $(\bm \Lambda_n)_{n\ge 1}$. Rewrite the allocation probability in (\ref{eqallocationPmul}) as follows
\begin{align}\label{eqallocationPnew}
 P(T_{n}=t|\bm{Z}_{n-1},Z_n=\bm k,\bm{T}_{n-1})
  = p_{r_t(\bm k)},
\end{align}
where $r_t(\bm k)=r\big(\Lambda_{n-1,t}(\bm k),\bm \Lambda_{n-1}(\bm k)\big)$  is the rank of $\Lambda_{n-1,t}(\bm k)$ among the elements of $\bm \Lambda_{n-1}(\bm k)=\big(\Lambda_{n-1,1}(\bm k), \ldots, \Lambda_{n-1,T}(\bm k)\big)$, i.e., $r_t(\bm k)=s$ if $\Lambda_{n-1,t}(\bm k)=\Lambda_{n-1,(s)}(\bm k)$. From (\ref{eqallocationPnew}), it follows that
\begin{equation}\label{eqallocationPnew2}\begin{matrix}
&\pr\Big(\Delta D_{n,t}(\bm k)=1-\frac{1}{T}\big|\mathscr{F}_{n-1}\Big)=p_{r_t(\bm k)}p(\bm k),\\
&\pr\Big(\Delta D_{n,t}(\bm k)=-\frac{1}{T}\big|\mathscr{F}_{n-1}\Big)=\big(1-p_{r_t(\bm k)}\big)p(\bm k),\\
&\pr\left(\Delta D_{n,t}(\bm k)=0 \big|\mathscr{F}_{n-1}\right)=1-p(\bm k).
 \end{matrix}
 \end{equation}
Then
\begin{align*} \ep[D_{n,t}(\bm k)|\mathscr{F}_{n-1}]= & D_{n-1,t}(\bm k) +\left[\Big(1-\frac{1}{T}\Big) p_{r_t(\bm k)}- \frac{1}{T} \big(1-p_{r_t(\bm k)}\big)\right]p(\bm k)\\
= & D_{n-1,t}(\bm k) +\left[ p_{r_t(\bm k)}- \frac{1}{T}\right]p(\bm k)\\
=& D_{n-1,t}(\bm k) + \overline{g}_{n-1,t,\bm k}p(\bm k),
\end{align*}
where $\overline{g}_{n-1,t,\bm k}=\overline{g}\big(\Lambda_{n-1,t}(\bm k),\bm \Lambda_{n-1}(\bm k)\big)=p_{r_t(\bm k)}- \frac{1}{T}$.
From (\ref{eqallocationPnew2}), it is easily seen that the transition probabilities of the Markov chain $(\bm D_n)_{n\ge 1}$ are symmetric about treatments, i.e., for any permutation $\Pi(1),\ldots,\Pi(T)$ of $1,\ldots, T$, the the transition probabilities of
$$
 \Big(D_{n,\Pi(1)}(\bm k),\ldots,D_{n,\Pi(T)}(\bm k)\Big)_{1\leq k_{1}\leq m_{1},\ldots, 1\leq k_{I}\leq m_{I}, n\ge 1}
$$
are the same. It follows that the  invariant probability measure  $\pi_{\lambda}$ of  $(\bm \Lambda_n)_{n\ge 1}$ is symmetric about treatments. So
$$ \pi_{\lambda}[p_{r_t(\bm k)}]=\frac{1}{T} \; \text{ and }  \pi_{\lambda}\overline{g}=0. $$
Now, by  Theorem 17.4.2 of Meyn and Tweedie (1993), there is a constant $R$ such that the Poisson's equation
\begin{equation}\label{PossionEquationMulti}
 \widehat{g}-P_{\lambda}\widehat{g}=\overline{g}-\pi_{\lambda}\overline{g}=\overline{g}
 \end{equation}
has a solution $\widehat{g}=\widehat{g}(t,\bm k)=\widehat{g}(t,\bm k,\bm\Lambda)$ which is a function of $\bm\Lambda$ defined on the state space of $\bm \Lambda$  with $|\widehat{g}|\le R(V+1)$, where $P_{\lambda}$ is the transition probability matrix of $\bm\Lambda$.
 Denote
$$
 \Delta M_{n,t,\bm k}=D_n(\bm k)-D_{n-1,t}(\bm k)+p(\bm k)[\widehat{g}(t,\bm k, \bm \Lambda_n)-\widehat{g}(t,\bm k,\bm \Lambda_{n-1})],
$$
 $M_{n,t}(i;k_i)=\sum_{\bm k\setminus k_i}M_{n,t,\bm k}$ and $g_{t,i;k_i}(\bm\Lambda_n)=\sum_{\bm k\setminus k_i}
 p(\bm k)[\widehat{g}(t,\bm k, \bm \Lambda_n)-\widehat{g}(t,\bm k,\bm \Lambda_0)]$.
 Then $\{M_{n,t,\bm k}\}$ and $\{M_{n,t}(i;k_i)\}$ are martingales, and
\begin{align}\label{eqDtoMartingale1}
&   M_{n,t,\bm k}= D_n(\bm k)-D_{0,t}(\bm k)+p(\bm k)[\widehat{g}(t,\bm k, \bm \Lambda_n)-\widehat{g}(t,\bm k,\bm \Lambda_0)]\\ \nonumber
&M_{n,t}(i;k_i)=D_{n,t}(i;k_i)-D_{0,t}(i;k_i)+g_{t,i;k_i}(\bm \Lambda_n).
\end{align}
With the same arguments as that in (\ref{eqVaRofM}) and (\ref{eqVaRofMMarg}) we can find constants $\sigma_t^2(\bm k)$ and $\sigma_t^2(i;k_i)$ such that
$$ \ep[M_{n,t,\bm k}^2]=n \sigma_t^2(\bm k)+O(1) \;\text{ and }\; \ep[M_{n,t}^2(i;k_i)]=n\sigma_t^2(i;k_i)+O(1). $$
It follows that
$$ \ep[D_{n,t}^2(\bm k)]=n \sigma_t^2(\bm k)+O(\sqrt{n}) \;\text{ and }\; \ep[D_{n,t}^2(i;k_i)]=n\sigma_t^2(i;k_i)+O(\sqrt{n}). $$

If $\sigma_t^2(\bm k)=0$, then with the same argument as in the two treatment case we   have $M_{n,t,\bm k}\equiv 0$. And then by (\ref{eqDtoMartingale1})
$$D_{n,t}(\bm k)-D_{0,t}(\bm k)=- p(\bm k)[\widehat{g}(t,\bm k, \bm \Lambda_n)-\widehat{g}(t,\bm k,\bm \Lambda_0)]$$
is a function of $\bm\Lambda_n$, which is a   contradiction when $w_s=0$, because one can find a path of $(\bm D_n)$ such that, at the end points of the path, the values of $(\bm \Lambda_n)$ are the same but the values of $\big(D_{n,t}(\bm k)\big)$ are different.

Similarly, if $\sigma_t^2(i;k_i)=0$, then we must have $M_{n,t}(i;k_i)\equiv 0$ and
$$ D_{n,t}(i;k_i)-D_{0,t}(i;k_i)=-g_{t,i;k_i}(\bm \Lambda_n), $$
which is also a    contradiction when $w_s+w_{m,i}=0$. The proof is now completed.
$\Box$
\vspace{2cm}


\bigskip
\baselineskip 15pt
\begin{center}
REFERENCES
\end{center}

\vspace{-0.1in} \footnotesize

\begin{enumerate}

\setcounter{enumi}{0}
\renewcommand{\theenumi}{{\rm [\arabic{enumi}]}}

\stepcounter{enumi}
\item[\theenumi]
Ashley, E. A., Butte,  A. J., Wheeler,  M. T., \textit{et al.}  (2010). Clinical evaluation incorporating a personal genome. \textit{The Lancet}  \textbf{375}, 1525-1535.


\stepcounter{enumi}
\item[\theenumi]
Baldi Antognini, A.  and Giovagnoli, A. (2004).
A new ¡®biased coin design¡¯ for the sequential
allocation of two treatments.  \textit{Journal of the Royal Statistical Society. Series C (Applied Statistics)}  {\bf 53}, 651-664.



\stepcounter{enumi}
\item[\theenumi]
Birkett, N. J. (1985). Adaptive allocation in randomized controlled trials. \textit{Controlled Clinical Trials} \textbf{6}, 146-155.

\stepcounter{enumi}
\item[\theenumi]
Ciolino, J., Zhao, W., Martin, R., \& Palesch, Y (2011). Quantifying the cost in power of ignoring continuous covariate imbalances in clinical trial randomization. \textit{Contemporary Clinical Trials} {\bf 32}, 250-259.

\stepcounter{enumi}
\item[\theenumi]
Efron, B. (1971). Forcing a sequential experiment to be balanced. \textit{Biometrika} \textbf{58}, 403-417.

\stepcounter{enumi}
\item[\theenumi]
Forsythe, A. B. (1987). Validity and power of tests when groups have been balanced for prognostic factors. \textit{Computational Statistics and Data Analysis} \textbf{5}, 193-200.



\stepcounter{enumi}
\item[\theenumi]
Hu, F. (2012). Statistical issues in trial design and personalized medicine. \textit{ Clinical Investigation} \textbf{2}, 121-124.

\stepcounter{enumi}
\item[\theenumi]
Hu, Y. and Hu, F. (2012).   Asymptotic properties of covariate-adaptive randomization. \textit{Annals of Statistics} {\bf 40}, 1794-1815.

\stepcounter{enumi}
\item[\theenumi]
Hu, F. and Rosenberger, W. F. (2006). {\em The Theory of Response-Adaptive Randomization in Clinical Trials}. John Wiley and Sons. Wiley Series in Probability and Statistics.

\stepcounter{enumi}
\item[\theenumi]
Hu, F. and Zhang L.-X. (2004). Asymptotic properties of doubly adaptive biased coin designs for
 multitreatment clinical trials. \textit{The Annals of Statistics} \textbf{32}, 268-301.

\stepcounter{enumi}
\item[\theenumi]
Hu, F., Zhang, L.-X. and He, X. (2009),  Efficient randomized adaptive designs. \textit{The Annals of Statistics} {\bf 37}, 2543-2560.



\stepcounter{enumi}
\item[\theenumi]
Kundt, G. (2009). Comparative evaluation of balancing properties of stratified
randomization procedures.  \textit{Methods of Information in Medicine} \textbf{48}, 129-134.

\stepcounter{enumi}
\item[\theenumi]
Lipkin, S. M., Chao, E. C., Moreno, V., \textit{et al.}  (2010).
Genetic variation in 3-hydroxy-3-methylglutaryl CoA reductase modifies the chemopreventive activity of statins for colorectal cancer.
\textit{Cancer Prevention Research}  {\bf 3}, 597-603.




\stepcounter{enumi}
\item[\theenumi]
Markaryan, T. and Rosenberger, W. F. (2010). Exact Properties of Efron's Biased Coin Randomization Procedure. \textit{Annals of Statistics} \textbf{38}, 1546-1567.


\stepcounter{enumi}
\item[\theenumi]
McIlroy M., McCartan D., Early, S., Gaora, P., Pennington, S. Hill, A. D. K. and Young, L. S. (2010). Interaction of developmental transcription factor HOXC11 with steroid receptor coactivator SRC-1 mediates resistance to endocrine therapy in breast cancer.
\textit{Cancer Research} \textbf{70}, 1585-1594.

\stepcounter{enumi}
\item[\theenumi]
Meyn, S. P.  and Tweedie, R. L. (1993). \textit{Markov chains and stochastic stability}.
Springer-Verlag, London.


\stepcounter{enumi}
\item[\theenumi]
Pocock, S. J. and Simon, R. (1975). Sequential treatment assignment
with balancing for prognostic factors in the controlled
clinical trial. \textit{Biometrics} \textbf{31}, 103-115.


\stepcounter{enumi}
\item[\theenumi]
Rosenberger, W. F. and Lachin, J. M. (2002). {\em
Randomization in Clinical Trials: Theory and Practice}. John Wiley and Sons. Wiley Series in Probability and Statistics.

\stepcounter{enumi}
\item[\theenumi]
Rosenberger, W. F. and Sverdlov, O. (2008). Handling covariates in the design of clinical trials. \textit{Statistical Science} \textbf{23}, 404-419.

\stepcounter{enumi}
\item[\theenumi]
Scott, N. W., McPherson, G. C., Ramsay, C. R. and Campbell, M. K. (2002). The method of minimization for allocation to clinical trials: A review. \textit{Control Clinical Trials} {\bf 23}, 662-674.

\stepcounter{enumi}
\item[\theenumi]
Shao,  J., Yu, X. and Zhong, B. (2010). A theory for testing hypotheses under covariate-adaptive
randomization. \textit{Biometrika} \textbf{97}, 347-360.

\stepcounter{enumi}
\item[\theenumi]
Signorini, D. F., Leung, O., Simes, R. J., Beller, E., Gebski, V. J. and Callaghan, T. (1993). Dynamic balanced randomization for clinical trials. \textit{Statistics in medicine} \textbf{12}, 2343-2350.



\stepcounter{enumi}
\item[\theenumi]
Taves, D. R. (1974). Minimization: A new method of assigning
patients to treatment and control groups. \textit{Clinical Pharmacology and Therapeutics} \textbf{15}, 443-453.

\stepcounter{enumi}
\item[\theenumi]
Taves, D. R. (2010). The use of minimization in clinical trials.
\textit{Contemporary Clinical Trials} \textbf{31}, 180-184.

\stepcounter{enumi}
\item[\theenumi]
Toorawa, R., Adena, M., Donovan, M., Jones, S. and Conlon, J. (2009). Use of simulation to compare the performance of minimization with stratified blocked randomization.  \textit{Pharmaceutical Statistics} \textbf{8}, 264-278.

\stepcounter{enumi}
\item[\theenumi]
Tymofyeyev, Y., Rosenberger, W.F. and Hu, F. (2007). Implementing
optimal allocation in sequential binary response experiments. \textit{
Journal of the American Statistical Association} \textbf{102}, 224-234.


\stepcounter{enumi}
\item[\theenumi]
Weir, C. J. and Lees, K. R. (2003). Comparison of stratification and adaptive methods for treatment allocation in an acute stroke clinical trial.  \textit{Statistics in Medicine}
\textbf{22}, 705-726.


\stepcounter{enumi}
\item[\theenumi]
Zhang, L.-X., Hu, F. and Cheung, S. H. (2006). Asymptotic theorems of sequential estimation-adjusted urn models for clinical trials.  \textit{The Annals of  Applied Probability} {\bf 16}, 340-369.

\stepcounter{enumi}
\item[\theenumi]
Zhang, L.X., Hu, F., Cheung. S.H. and Chan, W.S. (2007). Asymptotic
properties of  covariate-adjusted adaptive designs. \textit{Annals of
Statistics} {\bf 35}, 1166-1182.
\end{enumerate}

\end{document}